\documentclass[11pt]{amsart}
\usepackage{amssymb, amsmath, amsthm}
\usepackage[latin1]{inputenc}
\usepackage{graphicx}
\usepackage[final]{hyperref}
\usepackage{color}
\usepackage{ulem}
\usepackage{epstopdf}
\usepackage{cancel}

\usepackage[a4paper, centering]{geometry}
\usepackage{soul}
\usepackage{color}
\usepackage{graphicx}
\usepackage{scalerel,stackengine}

\usepackage{mathtools}

\newtheorem{theorem}{Theorem}
\newtheorem{proposition}[theorem]{Proposition}

\newtheorem{lemma}[theorem]{Lemma}

\theoremstyle{definition}

\theoremstyle{remark}
\newtheorem{remark}[theorem]{Remark}
\newtheorem{example}[theorem]{Example}
\definecolor{verde}{RGB}{20,150,100}
\definecolor{purple}{RGB}{200,30,200}

\newcommand{\EEE}{\color{black}}

\stackMath
\newcommand\reallywidecheck[1]{%
\savestack{\tmpbox}{\stretchto{%
  \scaleto{%
    \scalerel*[\widthof{\ensuremath{#1}}]{\kern-.6pt\bigwedge\kern-.6pt}%
    {\rule[-\textheight/2]{1ex}{\textheight}}
  }{\textheight}%
}{0.5ex}}%
\stackon[1pt]{#1}{\scalebox{-1}{\tmpbox}}%
}

\def\R{\mathbb{R}}

\newcommand{\sm}{\setminus}

\newcommand{\Om}{\Omega}
\newcommand{\sq}{\subseteq}
\newcommand{\ra}{\rightarrow}
\newcommand{\reflext}{{\mathcal R} _t} 
\newcommand{\reflexT}{{\mathcal R} _T}

\def \e{\varepsilon}

\begin{document}
\title[]{ 
Alexandrov theorem for general nonlocal curvatures: 
the geometric impact of the kernel
}

\author[]{Dorin Bucur, Ilaria Fragal\`a}

\thanks{}

\address[Dorin Bucur]{
Universit\'e  Savoie Mont Blanc, Laboratoire de Math\'ematiques CNRS UMR 5127 \\
  Campus Scientifique \\
73376 Le-Bourget-Du-Lac (France)
}
\email{dorin.bucur@univ-savoie.fr}

\address[Ilaria Fragal\`a]{
Dipartimento di Matematica \\ Politecnico  di Milano \\
Piazza Leonardo da Vinci, 32 \\
20133 Milano (Italy)
}
\email{ilaria.fragala@polimi.it}

\keywords{   Rigidity results,  nonlocal mean curvature, measurable sets, moving planes.  }
\subjclass[2010]{  53C24, 49Q15, 28A75.}
\date{\today}

\maketitle

\begin{abstract}  For a general radially symmetric, non-increasing, non-negative kernel $h\in L ^ 1 _{loc} ( \R ^ d)$, 
we study the rigidity of measurable sets in $\R ^ d$  with constant nonlocal $h$-mean curvature.
Under a suitable ``improved integrability'' assumption on $h$, 
we prove that these sets are finite unions of equal balls, 
as soon as they satisfy a natural nondegeneracy condition. Both the radius of the balls and their mutual distance can be controlled from below in terms
of suitable parameters depending explicitly on the measure of the level sets of $h$.  In the simplest, common case, in which $h$ is positive, bounded and decreasing, our result implies that any bounded open set or any bounded measurable set with finite perimeter which has constant nonlocal $h$-mean curvature has to be a ball.
\end{abstract}

\section{Introduction}\label{sec:main}

Let $h : \R ^ d \to \R _+$ be a radially symmetric  non-increasing measurable kernel. Given a measurable set  $\Omega  \subset \R ^d$, 
by seeing the quantity $h ( x-y)$ as an interaction density between two 
points $x \in \Omega$ and $y \in \Omega ^ c: = \R ^ d \setminus \Omega$,
 the  {\it nonlocal $h$-perimeter } of $\Omega$ can be defined by
$$P _ h (\Omega): = \int _\Omega \int _{\Omega^c} h (x-y) \, dx \, dy\,.$$

For typical choices of the kernel, the classical perimeter by Caccioppoli  (see \cite{Ma12}) can be recovered from the above nonlocal one by a scaling argument. 
Moreover,
by analogy with the classical case,  a natural notion of {\it nonlocal $h$-mean curvature} can be associated with the nonlocal $h$-perimeter, by setting
$$H _ h (\Omega):=  \int _{\R ^d} h ( x-y) \big ( \chi _{\Omega ^ c} (y) - \chi _\Omega ( y)) \, dy \,.$$ 
Actually, this definition is well-posed as soon as the regularity of $\Omega$ ensures the finiteness of the integral,  and 
is justified by the fact that minimizers of $P _ h$ under a volume constraint turn out to have constant $h$-mean curvature. 

The concept of nonlocal perimeter has been first considered in \cite{BBM01}, and it has been widely developed since then. 
In particular, based on the seminal papers \cite{CS08, CRS10}, a wide attention has been devoted to the so-called {\it fractional} perimeter, 
which corresponds to the choice of the singular kernel
\begin{equation}\label{f:frach} h( x) = \frac{1}{|x| ^ {d+s} }  \, \qquad s \in (0, 1)\,. 
\end{equation}
Research  in the fractional setting has been extended to a
broad spectrum of directions, including
the study of isoperimetric type inequalities \cite{FS08}, minimal surfaces \cite{DSV17, FV17}, diffusion processes \cite{BPSV14, BV16, BSV17},
and mean curvature flows \cite{CMP15, CNR17, I09} (where references are clearly a sparse sampling). 

These topics have been investigated also for another family of kernels, namely the one of bounded integrable kernels, 
 which has been extensively treated in  the monograph \cite{MRT}. A  distinguishing feature of  such kernels, that 
we care to mark as a deep difference from the fractional one, is that the notion of nonlocal mean curvature makes sense for general  measurable sets. 

The present work is focused on  a recent trend in nonlocal analysis, namely the rigidity of sets with constant nonlocal mean curvature. 
The reference milestone result in the local setting is  Alexandrov theorem, dating back to 1958: it states that, among  connected smooth domains, the only one with constant mean curvature is the ball \cite{ale}.  Let us also mention that, still in the local framework, 
a significant extension of Alexandrov theorem has been obtained in the recent paper \cite{DelMag} by
Delgadino and Maggi, who have been able to remove any  kind of regularity or connectedness assumption: 
they have proved that,  among sets with finite measure and finite perimeter, the only ones with constant mean curvature (now meant in distributional sense) are finite unions of equal balls.   

In the nonlocal setting,  the problem has been attacked in recent years for two distinct choices of the kernel, that we shortly summarize hereafter. 

The first case, which has been solved in 2018, is that  of the  fractional kernel \eqref{f:frach}: 
in two independent papers, Ciraolo-Figalli-Maggi-Novaga \cite{CFMN18} and Cabr\'e-Fall-Sol\`a Morales-Weth  \cite{CMMW} have proved that, among sets of class  $C ^ { 1, \alpha} $, the only one with constant fractional mean curvature is the ball. 
Notice that the regularity assumption in such result cannot be removed, since it is necessary to give a meaning to the fractional mean curvature; on the other hand, no connectedness hypothesis is needed,  because nonlocal interactions automatically rule out the bubbling phenomenon appearing in \cite{DelMag}. 

The second case, which has been treated in our previous work \cite{BF},  is that of a completely different kernel, 
given by
\begin{equation}\label{f:ballh}
h = \chi _{ B _ r (0)}\,,
\end{equation} 
where $r$ is a {\it fixed} positive radius, and $B _ r (0)$ denotes the ball of radius $r$ centred at the origin. 
Since such kernel is nonsingular, the main novelty of the corresponding rigidity result is that its validity extends to the broad class of measurable sets.  Still as a consequence of the kernel's properties, specifically of the boundedness of its support,  some kind of ``short-range connectedness'' assumption, that we call  $r$-{\it nondegeneracy}, is needed in order to get rigidity.  Under this additional assumption (which holds for free e.g.\ for open connected sets of diameter larger than $r$), we proved that the only measurable sets with constant nonlocal mean curvature for the kernel \eqref{f:ballh}  are finite unions of equal balls, of radius  $R> r/2$, lying at distance at least $r$ from each other. 
Thus the initial choice of the positive radius $r$  tunes the rigidity phenomenon, in the sense that  $r$ acts as a threshold from below both for the diameter of the balls, and for their mutual distance.

Aim of this paper is to go beyond in the study of rigidity  in the nonlocal setting, 
by   
characterizing  measurable sets with constant $h$-mean curvature  when $h$ is an {\it arbitrary} locally integrable kernel.     
These sets, that we call   {\it $h$-critical}, satisfy  the following condition: 
\begin{equation} \label{f:hyp3}
\exists   c >0 \ :\ \int_\Omega h  (x-y) \, dy  = c \qquad \forall x \in {\partial^* \Omega}\,,
\end{equation} 
where $\partial ^* \Omega$ denotes the essential boundary of $\Omega$, namely the set of  points $x\in \R^d$  at which  both $\Omega$ and its complement $\Omega ^ c$ have a strictly positive $d$-dimensional upper density.  

Apart from the local integrability assumption on $h$, which is
the minimal requirement to make the $h$-mean curvature well-defined for all measurable sets, 
our scope is to analyse the problem in a unique framework, including  both  singular and nonsingular kernels; 
in particular, the challenge is to throw some light on the delicate interplay between the geometry of the kernel and the corresponding rigidity result. This aspect  requires to go into depth in the comprehension of the principles which govern the rigidity phenomenon. 

A crucial issue is that, 
 in order to get rigidity, the $h$-criticality condition must be combined with a sort of ``short-range connectedness'', which extends in a natural way the above mentioned notion of $r$-nondegeneracy introduced in \cite{BF}. Precisely, we say that a measurable set is $h$-{\it nondegenerate} if 
\begin{equation}\label{f:hyp4} \inf _{x_1, x_2 \in \partial ^* \Omega} 
\frac{\int_{\Omega}  |h ( x_1 - y) - h ( x _ 2 -y) |      }{\| x_1 -x_2 \| }  \, dy  > 0\,. 
\end{equation}

Some sufficient conditions for nondegeneracy  will be given in Proposition \ref{p:cake}. For instance, if $\Omega$ is an open connected set, or an indecomposable set of
 finite perimeter which is $h$-critical, condition \eqref{f:hyp4} is  fulfilled as soon as the diameter of $\Omega$ is larger than the radius of the  the level set $\{h ={\rm ess}\, {\rm sup}\,  h\}$: this means in particular that, if $h$ does not have a plateau of positive measure at its supremum, 
such a set $\Omega$ is automatically $h$-nondegenerate.

In order to get rigidity  for measurable sets which are $h$-critical and $h$-nondegenerate, we also need to ask  a technical condition on $h$, which is an
{\it improved integrability assumption}, expressed through the
behaviour of its level sets. We ask that 
\begin{equation}\label{f:integrability} 
  \int _ {1}  ^ { + \infty}\!\!\!   r ^ { d-1  }(s)  \, ds  < + \infty \,,
\end{equation}
where $r ( s)$ is the distribution function of the map  $\varphi : \R _+ \to \R _+$  defined by
\begin{equation}\label{f:varphi}
h ( x) = \varphi ( |x|) \qquad \forall x >0\,, 
\end{equation}  
namely
\begin{equation}\label{f:r}
r ( s):= \mathcal L ^ 1  \big ( \big \{ x \in \R_+  \ :\  \varphi  ( x) > s \big \} 
\big)  \qquad \forall s \geq 0 \,.
\end{equation}

Condition \eqref{f:integrability} is satisfied in particular when $h$ is bounded,  or when it  has a ``not too steep'' singularity at $0$ (see the examples at the end of this Introduction). 

\smallskip
 
Our rigidity result reads: 

\begin{theorem} \label{t:serrin3} Let $h$ be a radially symmetric, non-increasing, non-negative kernel in  $L ^ 1 _{\rm loc} (\R ^d)$ satisfying the improved integrability condition \eqref{f:integrability}. 
Let $\Omega$ be a set of finite Lebesgue measure which is $h$-critical and $h$-nondegenerate, i.e.\ satisfies \eqref{f:hyp3} and \eqref{f:hyp4}. Then 
$\Omega$ is equivalent to  
a finite union of balls $B _ i$ of the same radius $R$; moreover,   
\begin{eqnarray}
& \displaystyle R> \frac{\eta}{2} \,, \quad \text{ with } \   \eta :=  \mathcal L ^ 1 \big  ( \{ \varphi =  {\rm ess}\, {\rm sup} \, \varphi \}  \big ) 
\,,\qquad \quad \quad
& \label{f:raggio}
\\ \noalign{\bigskip}
& {\rm dist} (B _i, B _ j)  \geq  r(\sigma) \, , \quad\text{ with }  \sigma : = \begin{cases}
 0 & \text{  if }    {\rm diam}\, \Omega \geq r(0) 
\\  \noalign{\medskip}
\varphi ({\rm diam} \Omega) & \text{  if }  {\rm diam}\, \Omega <  r(0)  \,.
 \end{cases}
& \label{f:distanza}
\end{eqnarray} 

  \end{theorem}

\begin{remark}
Let us point out that, in the simplest, common case, in which $h$ is positive, bounded and decreasing, Theorem  \ref{t:serrin3} implies that any bounded open set or any bounded measurable set with finite perimeter which has constant nonlocal $h$-mean curvature has to be a ball. In more general situations, the interpretation of Theorem \ref{t:serrin3} (especially concerning the role of the two parameters $\eta$ and $\sigma$) 
is discussed below and in the next sections.
 \end{remark} 

\begin{remark}[on the size of balls]  The parameter $\eta$  depends just on the kernel:
in particular, the positivity of $\eta$ occurs only when  the kernel is bounded and reaches its supremum on a
plateau of positive measure. Thus we can distinguish two cases:

\smallskip 
--  Case $\eta= 0 $  (no supremal plateau): balls may have  arbitrarily small scale.

\smallskip
--  Case $\eta>0 $  (supremal plateau): the diameter of balls is bounded from below by  the radius  of the plateau. 

 \end{remark}

\begin{remark}[on the mutual distance of multiple balls]  The parameter $\sigma$ depends on the interplay between the diameter of $\Omega$ and the radii of the level sets of the kernel.  More precisely,  after noticing that $r(0) = \mathcal L ^ 1 {(\rm supp} \, \varphi)$, we can distinguish two cases: 

\smallskip
-- Case ${\rm diam}\, \Omega \geq \mathcal L ^ 1 {(\rm supp} \, \varphi) $ ($\Omega$ is  ``large'' compared to the support of  the kernel): multiple balls  are allowed, at mutual distance bounded from below by $ \mathcal L ^ 1 {(\rm supp} \, \varphi)$  
(that is, by
the radius of ${\rm supp} \, h$).

\smallskip
-- Case 
${\rm diam}\, \Omega <  \mathcal L ^ 1 {(\rm supp} \, \varphi) $  ($\Omega$ is ``small'' compared to the support of the kernel): by the properties of the distribution function (see Lemma \ref{l:distribution})  it holds
$$
\mathcal L ^ 1  \big (\{ \varphi >  \varphi ( {\rm diam} \Omega )  \}   \big ) = r (\varphi ({\rm diam} \, \Omega) ) \leq {\rm diam } \Omega
\leq  r( \varphi({\rm diam} \Omega) ^ - ) = \mathcal L ^ 1  \big ( \{ \varphi  \geq \varphi ( {\rm diam} \Omega )   \}   \big ).$$
Hence,  a necessary condition for $\Omega$ to be a 
{\it multiple} family of balls   is that the kernel has  a  level set  of positive measure, corresponding to a jump in the distribution function $r$. 
Two subcases may occur:

\begin{itemize}\item[--]  Case $\mathcal L ^ 1  \big ( \{ \varphi  = \varphi ( {\rm diam} \Omega )   \}   \big ) = 0$:  
multiple balls are not allowed, because by inequality \eqref{f:distanza} they should be at  distance equal at least to  ${\rm diam}\,  \Omega$. 

\smallskip
\item[--] Case $\mathcal L ^ 1  \big ( \{ \varphi  = \varphi ( {\rm diam} \Omega )   \}   \big ) > 0$: 
multiple balls are allowed, at mutual distance bounded from below by $r (\varphi ({\rm diam} \Omega) )$, and {\it all of them will be contained into the same level set of $h$}, given by  points $x \in \R ^d$ such that $
\varphi(|x|) \geq \varphi ({\rm diam} \Omega)$.

\end{itemize} 
\end{remark}

\begin{example}
Let
$$h ( x)=  \frac{1}{|x| ^{\alpha}}  \qquad  \text{ with } \ \alpha < d -1 \,.$$
Due to the choice of $\alpha$, the improved integrability condition \eqref{f:integrability} is satisfied. 
Theorem \ref{t:serrin3} applies: we have $\eta= 0$,  and $h$ does not have any level set of positive measure; 
hence, the unique $h$-critical and $h$-nondegenerate domain of finite measure is a single ball, 
whose radius can be arbitrarily small. \end{example}

\begin{example}
Let 
$$h ( x)= \sum _ {i = 1} ^N \alpha_i \chi _{B _ {r_i} (0)  \setminus B _ {r _ {i-1} } (0)  }\,.
$$
where $\alpha _ 1 > \alpha _ 2 > \dots > \alpha _N >0$ and $0= r_0  < r _ 1 < r _2 < \dots < r _N$. 
Since $h$ is bounded, the improved integrability condition \eqref{f:integrability} is satisfied. By
Theorem \ref{t:serrin3}, a domain of finite measure which is $h$-critical and $h$-nondegenerate 
will be a finite union of equal balls, of  radius $R> 
\eta= r_1/2$, with 
$${\rm dist} (B_i, B _j) \geq \sigma =   \begin{cases}
  \noalign{\medskip} r_N & \text{  if }    {\rm diam}\, \Omega \geq r_N
\\ \noalign{\medskip}
r_{i}  & \text{  if } r_{i} \leq  {\rm diam}\, \Omega < r_{i+i}  \text{ for } i = 1, \dots, N \,.
 \end{cases}
$$
In particular, for $N = 1$, we recover the result proved in \cite{BF}. 
\end{example} 

 The proof of Theorem \ref{t:serrin3} is obtained 
via a new version of the moving planes method
valid in the framework of measurable sets,  which has been settled in \cite{BF}.
However, with respect to the case $h = \chi _{B _ r (0)}$, 
dealing with  a nonconstant kernel makes the proof considerably more  delicate. 
This is the reason why we decided to omit all the parts of the proof which closely follow \cite{BF}, 
and in spite to focus in full detail on all the parts where the kernel plays an important role. 
Such parts are attacked relying on the basic idea of layering integrals according to Cavalieri's principle:  in particular, this allows to set up some key estimates, which require the improved integrability assumption \eqref{f:integrability}.

The paper is organized as follows: the required preliminaries are collected in Section \ref{sec:prel},  and then the proof is given in Section \ref{sec:proof} (which is  in turn divided, for the sake of clearness, into four subsections).

\section {Preliminaries}\label{sec:prel} 

Throughout the paper,  we assume that $h$ is a radially symmetric, non-increasing, non-negative measurable function in  $L ^ 1 _{\rm loc} (\R ^d)$. Moreover, for any $x \in \R ^ d$, we set for brevity 
$$h _ x (y ):= h ( x-y) \qquad \forall y \in \R ^ d\,.$$ 

 \subsection {On some plain consequences of criticality}\label{sec:plain}

    \begin{lemma}\label{l:bound}   
    Let  $\Omega$ be a measurable set of finite Lebesgue measure satisfying the criticality condition \eqref{f:hyp3}. 
Then $\Omega$ is bounded. 
  \end{lemma} 
    
    \proof Assume without loss of generality that $\inf _ \R h = 0 $ (otherwise replace $h$ by $h - \inf _\R h$).  By contradiction, let $\{ p _n \} $ be a sequence of points in $\partial ^ * \Omega$, with $|p _n| \to + \infty$.  Since $|\Omega | < + \infty$, for every $\varepsilon >0$ there exists $R_\varepsilon$, with $R_\varepsilon \to + \infty$ as $\varepsilon \to 0$,  such that $|\Omega \cap B _ {R_\varepsilon} (p_n ) |<  \omega _ d \varepsilon ^ d $. Thus we have
    $$\int _\Omega h _ { p _n} =  \int _{\Omega  \setminus B _ {R_\varepsilon} ( p _n )} h _ { p _n}  + \int _{\Omega \cap B _ {R_\varepsilon} ( p _ n)} h _ { p _n}  \leq h ( R_\varepsilon ) |\Omega  | + \int _{B _ \varepsilon ( 0 ) } h  \,. $$  
 In the limit as $\varepsilon \to 0$, since $h \in L ^ 1 _{\rm loc} (\R ^d)$, this contradicts assumption \eqref{f:hyp3}. 
\qed

\bigskip

\begin{lemma}\label{l:chiusura} 
    Let $\Omega$ be a measurable set of finite Lebesgue measure. 
    \begin{itemize}
    \item[--] If  $\Omega$ satisfies the criticality condition \eqref{f:hyp3}, the same equality continues to hold 
    at every point $x \in \overline { \partial ^* \Omega}$. 
    
    \smallskip
\item[--]   If  $\Omega$ satisfies the nondegeneracy condition \eqref{f:hyp4}, the same strict inequality continues to hold
when the infimum is taken over the pairs  of points $x_1, x_2 \in \overline { \partial ^* \Omega}$.   
\end{itemize}
\end{lemma} 
\proof  Assume that $\Omega$ satisfies the criticality condition \eqref{f:hyp3}.
Let  $x_0= \lim _n x _n$, with $x_n \in \partial ^* \Omega$. Let us prove that \eqref{f:hyp3} continues to hold at $x = x_0$.  Set $\Omega _n = x_n - \Omega$ and $\Omega _0 = x_0 - \Omega$.  
We claim that, up to passing to a (not relabeled) subsequence,  
\begin{equation}\label{f:pointwise}
\chi _{\Omega _n} \to \chi _{\Omega _0} \qquad \text{pointwise  a.e. in } \R ^ d\,.
\end{equation}  
Indeed, since by assumption  $\Omega$ has finite Lebesgue measure, the sequence $\chi _{\Omega _n}$ is  bounded  in $L ^ 2 ( \R ^d)$ and hence, up to a subsequence, it converges weakly  in $L ^ 2 ( \R^d)$ to some function $f$. 
Since, for every $\varphi \in   \mathcal C ^ \infty _0 ( \R ^ d)$, 
$$\int _{\R ^ d } \chi_{\Omega _n } (y)  \varphi (y) \, dy  = 
\int _{x_n - \Omega }   \varphi (y) \, dy  =\int _{\Omega }   \varphi (x_n-y) \, dy =
\int _{\R ^ d } \chi_{\Omega  } (y) \varphi  (x_n - y )\, dy ,$$  
by dominated convergence we infer that the weak limit $f$ agrees with $\chi _{\Omega _0}$. 
Taking into account that
$$\| \chi _{\Omega _n } \| _{ L ^ 2 (\R ^ d)  } = |\Omega _n | = |\Omega _0| = \| \chi _{\Omega _0 } \| _{ L ^ 2 (\R ^ d)  }  \,,$$
we deduce that the convergence is strong in $L ^ 2 ( \R ^ d)$. 
Up to choosing a further subsequence, we have that \eqref{f:pointwise} is satisfied. 

Now consider the sequence $h _n := h \chi _{\Omega _n} $. By \eqref{f:pointwise}, up to  a subsequence it converges to  $h _0:= h \chi _{\Omega _0}$ pointwise a.e.\ in $\R ^ d$.  Moreover, we have  $|h _ n| \leq  |h|$. Since by assumption $h \in L ^ 1_{\rm loc} (\R ^ d)$,  and since by Lemma \ref{l:bound} $\Omega$ is bounded, by applying the dominated convergence theorem on a sufficiently large ball, we infer that 
$h _ n \to h_0$ in $L ^ 1 ( \R ^ d)$.  Hence, 
$$ \|h _0 \| _{L^ 1 (\R ^ d) }  = \lim _ n  \|h _n \| _{L^ 1 (\R ^ d) }  = c \,, $$
as required.  The proof of the second claim in the statement is analogous.   \qed

\subsection {On the distribution function and layered integrals}\label{sec:layer} 

For convenience of the reader, we recall in the next lemma some well-known properties of the function $r$ which maps any $s \in \R _+$ into the radius of the level set $\{ h > s\}$ (see for instance \cite{CF02} and references therein). 

\begin{lemma}\label{l:distribution}  Let $\varphi$ be associated with $h$ as in \eqref{f:varphi}, and $r ( s) $  be its distribution function defined in \eqref{f:r}.   Then: 
 
 \smallskip
 (i) The map $s \mapsto r ( s) $ 
   is non-increasing, with
$r ( s) = 0$ for every  $s \geq {\rm ess}\, {\rm sup} \, \varphi$.

\smallskip
(ii) The map $s \mapsto r ( s) $ it is right continuous and, setting $ r ( s_0^-):= \lim _{ s \to s_0^-} r ( s)$, it holds 
$$ r ( s_0^-) - r (s_0) =   \mathcal L ^ 1 (  \{ \varphi  \geq s_0 \} )  -  \mathcal L ^ 1 (  \{ \varphi >  s_0 \} )   =\mathcal L ^ 1(  \{ \varphi = s_0 \} )   \,.$$

in particular, $r$ it  is continuous at a given point $s_0$ if and only if $\mathcal L ^ 1(  \{ \varphi = s_0 \} )  = 0 $; 

(iii) It holds
$$ \sup _{(0, + \infty)} r (s) = r (0 ) = \mathcal L ^ 1 (\{ {\rm supp} \, \varphi\}) \quad \text{ and } \quad 
\inf _{(0, {\rm ess}\, {\rm sup} \, \varphi )} r (s) = \eta:= \mathcal L ^ 1 \big  ( \{ \varphi =  {\rm ess}\, {\rm sup} \, \varphi \}  \big ) 
 \,.$$ 
 
 (iv) It holds
 $$\varphi ( t) = \sup \big  \{ s \geq 0 \ :\  r ( s) > t \big   \} \qquad \forall t \geq 0 \,.$$ 
\end{lemma} 

\medskip
\begin{remark}
We point out for later use that  the parameters $\eta$ and $\sigma$ introduced in the statement of Theorem \ref{t:serrin3} enjoy the following properties, which can be easily checked by  using the above lemma:
 \begin{eqnarray} \forall \lambda >0 \, , \quad \mathcal L ^ 1 \big ( \big \{ s \, :\, r ( s)  \in (0, \lambda)   \big \}  \big ) > 0 \  \Leftrightarrow \  
\lambda >  \eta   & \label{f:misura0}  
\\ \noalign{\medskip}
r ( s) > {\rm diam} \, \Omega \  \ \Leftrightarrow \ \ s < \sigma \,.  \qquad \qquad \qquad & \label{f:sigma}
 \end{eqnarray}
\end{remark}

\bigskip
 In the following simple lemma, which will be  used repeatedly in the sequel, we exploit the layer-cake principle to rewrite the integrals appearing in the criticality and in the nondegeneracy conditions  in terms of the function $r ( s)$ and of the parameter $\sigma$.

\begin{lemma}\label{l:cake} Let $\Omega$ be a measurable set of finite Lebesgue measure satisfying the criticality condition \eqref{f:hyp3}, and let $\sigma$ be defined as in \eqref{f:distanza}. 

\smallskip
(i)  For every  $x  \in \overline{\partial ^* \Omega}$, it holds
\begin{equation} \label{f:cake1}  \int _\Omega h _ x ( y) \, dy  
  =  \int _0 ^ {+ \infty} |\Omega \cap  {B _ { r (s) } (x)}  | \, ds 
  = \sigma |\Omega| + \int _{\sigma} ^ {+ \infty} |\Omega \cap  {B _ { r (s) } (x)}  | \, ds  
 \,.
   \end{equation}
   
   \smallskip 
(ii)  For every $x _1, x _ 2 \in \overline{{\partial ^* \Omega}}$,  it holds 
\begin{equation}\label{f:cake2} 
 \int _\Omega  | h _ { x_1} - h _ { x_2}| \, dy =\int _{\sigma} ^ {+ \infty}  
\big | \Omega  \cap
\big ( B _ { r (s) } (x_1) \Delta B _ { r (s) } (x_2) \big )  
 \big   | \, ds  \,.
 \end{equation} 

\end{lemma} 

\proof 

For every $x \in \overline{\partial ^* \Omega}$, since by assumption $h _ x \in L ^ 1 _{\rm loc} (\R ^d)$, and since by Lemma \ref{l:bound} the set $\Omega$ is bounded, we have $h_x \in L ^ 1 (\Omega)$. Then,   by the layer-cake principle and Fubini Theorem,  we have
$$\begin{array}{ll} \displaystyle \int _\Omega h _ x ( y) \, dy  & \displaystyle = \int _\Omega \int _0  ^ {+ \infty} \!\!\! \chi _ {\{ h _ x > s \} } \, ds \, dy 
 =  \int _0  ^ {+ \infty} \!\!\! \int _\Omega \chi _ {\{ h _ x > s \} } \, dy \, ds
 \\ \noalign{\medskip} 
 & \displaystyle =  \int _0  ^ {+ \infty}\!\!\!\! \int _\Omega \chi _ {B _ { r (s) } (x)} (y) \, dy  \, ds =  \int _0  ^ {+ \infty} |\Omega \cap  {B _ { r (s) } (x)}  | \, ds  \  \,.\end{array} 
$$ 
The equality \eqref{f:cake1} follows by noticing that, as a consequence of \eqref{f:sigma}, for $s < \sigma$ we have  $|\Omega \cap  {B _ { r (s) } (x)}  | = |\Omega|$.

\smallskip 
For every  $x _1, x _ 2 \in \partial ^* \Omega$, we  
have 
$$\int _\Omega  | h _ { x_1} - h _ { x_2}| \, dy = \int _{\Omega  \cap \{ h _ { x_1} > h _ {x_2} \}}    ( h _ { x_1} - h _ { x_2})  \, dy  + 
 \int _{\Omega  \cap \{ h _ { x_2}  >  h _ {x_1} \}}   ( h _ { x_2} - h _ { x_1})  \, dy  =: I '  + I'' \,.$$ 
 
As above, we  use the layer cake principle to rewrite  $I ' $ as 
$$\begin{array}{ll}  I ' & \displaystyle = \int _0  ^ {+ \infty} \!\!\!  \int _{\Omega  \cap \{ h _ { x_1} > h _ {x_2} \}}
\big [ \chi _ {B _ { r (s) } (x_1) }    (y) - \chi _ {B _ { r (s)}  (x_2)  }  (y)    \big ] \, ds  
\\  \noalign{\medskip} 
& \displaystyle =\int _0 ^ {+ \infty}  \big |\Omega  \cap \{ h _ { x_1} > h _ {x_2} \} \cap
\big ( B _ { r (s) } (x_1) \setminus B _ { r (s) } (x_2) \big )  
 \big | \, ds \\ \noalign{\medskip} 
& \displaystyle  - \int _0  ^ {+ \infty}  \big | \Omega  \cap \{ h _ { x_1} > h _ {x_2} \} \cap \big ( B _ { r (s) } (x_2) \setminus B _ { r (s) } (x_1) \big )   \big | \, d s   \,.
\end{array} 
$$ 
Then, since $h$ is non-increasing, and since $\Omega \subset B _ { r (s) } (x_2)$ for $s > \sigma$ (again by \eqref{f:sigma}), 
we obtain 
$$\begin{array}{ll}  I ' & \displaystyle = \int _0  ^ {+ \infty}  \big | \Omega  \cap
\big ( B _ { r (s) } (x_1) \setminus B _ { r (s) } (x_2) \big )  
 \big   | \, ds 
 \\ \noalign{\medskip} 
& =  \displaystyle\int _{\sigma}   ^ {+ \infty}  \big | \Omega  \cap
\big ( B _ { r (s) } (x_1) \setminus B _ { r (s) } (x_2) \big )  
 \big   | \, ds 
  \,.
\end{array}
$$ 
Likewise, we obtain 
$$I'' = \int _{\sigma}  ^ {+ \infty}   \big | \Omega  \cap
\big ( B _ { r (s) } (x_2) \setminus B _ { r (s) } (x_1) \big )  
 \big   | \, ds \,. $$ 
 By adding the above expressions for $I'$ and $I''$, we obtain \eqref{f:cake2}. 
 \qed 
 
 \bigskip 

\bigskip

\subsection {On some sufficient conditions for nondegeneracy }

\begin{proposition}\label{p:cake}
Let $\Omega$ be a measurable set of finite Lebesgue measure satisfying the criticality condition \eqref{f:hyp3}. 
Assume in addition that $\Omega$ is either open or of finite perimeter. Then the nondegeneracy condition \eqref{f:hyp4}  is satisfied provided
 \begin{equation}\label{f:controimmagine} 
\inf _i {\rm diam} (\Omega_i)  > \eta \,, 
\end{equation} 
where $\{\Omega_i\} _i$  the family  of the connected or indecomposable components of $\Omega$, and  $\eta$ is defined as in \eqref{f:raggio}\,.
 \end{proposition}

\proof Working component by component, we are reduced to show that, if $\Omega$ is  an indecomposable set of finite perimeter, or an open connected set, it is not $h$-degenerate provided
\begin{equation}\label{f:controimmagine_bis}      {\rm diam} (\Omega)  > \eta \,.
\end{equation}  
By Lemma \ref{l:cake}, we have
$$ \inf _{x_1, x_2 \in \partial ^* \Omega} 
 \frac{\int_{\Omega}  \big |h ( x_1 - y) - h ( x _ 2 -y) |   \big |  }{\| x_1 -x_2 \| }  \, dy  =  
 \inf _{x_1, x_2 \in \partial ^* \Omega}  \int _{\sigma}  ^ {+ \infty}   \frac{\big | \Omega  \cap
\big ( B _ { r (s) } (x_1) \Delta B _ { r (s) } (x_2) \big )  
 \big   |} {\| x_1 -x_2 \| } \, ds\,.$$ 
 Passing the infimum under the sign of integral, we get
$$ \inf _{x_1, x_2 \in \partial ^* \Omega} 
\frac{\int_{\Omega} \big |h ( x_1 - y) - h ( x _ 2 -y) |   \big |  }{\| x_1 -x_2 \| }  \, dy  \geq \int _{\sigma}  ^ {+ \infty}  
\!\!\! \!\!\!\inf _{x_1, x_2 \in \partial ^* \Omega}   \frac{\big | \Omega  \cap
\big ( B _ { r (s) } (x_1) \Delta B _ { r (s) } (x_2) \big )  
 \big   |} {\| x_1 -x_2 \| } \, ds  \,.$$
 The
r.h.s.\ of the above inequality is strictly positive provided 
\begin{equation}\label{f:misura} \mathcal L ^ 1 \big ( \big \{ s \ :\ r ( s) \in ( 0,   {\rm diam} (\Omega)) \big \} \big ) >0 \,,
\end{equation}
because, by \cite[Proposition 10]{BF}, 
$$ \inf _{x_1, x_2 \in \partial ^* \Omega}   \frac{\big | \Omega  \cap
\big ( B _ { r } (x_1) \Delta B _ { r } (x_2) \big )  
 \big   |} {\| x_1 -x_2 \| }>0 \qquad \forall r \in (0, {\rm diam} (\Omega))\,.$$
In turn,  recalling \eqref{f:misura0},  condition \eqref{f:misura} is fulfilled thanks to assumption \eqref{f:controimmagine_bis}. 
   \qed

  \section{ Proof of Theorem \ref{t:serrin3}} \label{sec:proof} 
  
  {\bf Outline.}  We adopt the moving planes method for measurable sets settled in \cite{BF}. 
  
  As in the classical moving planes method,  the idea is to consider,  
for any fixed direction $\nu \in S ^ {d-1}$,  an initial hyperplane $H _ 0$ with unit normal $\nu$, not intersecting $\overline {\partial ^* \Omega}$ (this can be done thanks to Lemma \ref{l:bound}).   Then one starts moving $H _ 0$ in the direction of its normal $\nu$ to new positions $H _t$, 
so that at a certain moment of the process it starts intersecting $\overline {\partial ^* \Omega}$. 
  The main novelty of the approach introduced in \cite{BF} with respect to the classical case is how to define the stopping time of the movement, and then how to get rigidity (including possibly multiple balls). 

Here we follow  the same global strategy as  in \cite{BF}, but each part of the proof needs to be significantly changed,  due to the much greater generality of the kernel we work with (in the few points where  the same arguments apply, 
the reader is explicitly referred to \cite{BF}).

Before starting,  let us fix some notation and terminology: we denote by $H _ t^-$ and $H _ t ^ + $ the two closed halfspaces determined by $H _t$ (for definiteness, assume that $H _0 \subset H _ t ^-$); we  set 
 $$\Omega_ t  := \Omega \cap H _ t ^ {-}
 \, , \qquad \reflext := \text { the reflection of $\Omega _t $ about $H _t$}. 
$$

\smallskip

$\bullet$ We say that {\it symmetric inclusion} holds at $t$ if
 \begin{equation}\label{f:moving} 
\reflext \subset \Omega \qquad \text{ and } \qquad \Omega _t   \cup \reflext   \text{ is Steiner symmetric about $H _t$\,.} \end{equation} 
(Recall that a  measurable set   $\omega$   is Steiner symmetric about a hyperplane $H$ with unit normal $\nu$ if 
it is equivalent to the set of points
$x \in \R ^d$ of the form $x = z + t \nu$, with $z \in H$ and $|t| < \frac{ 1}{2} \mathcal H ^ 1 \big (
 \omega \cap \big \{ z + t \nu  : t \in \R \big \} \big )$).    

\smallskip 
$\bullet$ We say that  symmetric inclusion occurs at $t$ if  {\it with away contact} if \eqref{f:moving} holds and there exists an ``away contact point'', namely a point
 \begin{equation}\label{f:touching}
 p' \in  \big [ \overline {\partial ^* \reflext } \cap \overline {\partial ^*  \Omega}\big ] \setminus H _ t  \, .
  \end{equation}
when \eqref{f:moving} holds but \eqref{f:touching} is false, we say that symmetric inclusion at $t$ holds  {\it without away contact}. 

\smallskip
$\bullet$ We say that symmetric inclusion occurs at $t$ {\it with close contact}  if \eqref{f:moving} holds and there exists a ``close contact point'', namely a point  \begin{equation}\label{f:touching2}
H _ t \ni  q = \lim _n q _{1,n} = \lim _n q  _ {2,n},   \quad  q _{i, n}  \in  \overline{\partial ^* \Omega}\cap \{ q + t \nu\, :\, t \in \R \}, \quad   q _ {1,n} \neq  q _{2,n}\,.  \end{equation}
  Notice that symmetric inclusion can occur at the same $t$ with both away contact and close contact.

 \bigskip

We are now ready to start. We proceed in four steps, which  are carried over in separate subsections below.  

\subsection{Step 1 (start).}   
We prove the following claim: 

\smallskip
 $\bullet$ {\it Claim 1: There exists $\e >0$ such that, for every $t \in [0, \e)$, symmetric inclusion holds.  }

\smallskip 
The proof is based on the following  \EEE 

\begin{lemma}[no converging pairs]\label{l:nopairs} 
Under  the assumptions of Theorem \ref{t:serrin3}, 
if $\Omega $ is contained  into
$ H _0 ^+:=\{ z + t \nu\, :\, z \in H _0\, , \ t \geq 0 \}$, 
 $H_0$ being a hyperplane  with unit normal $\nu$,  there cannot exist two sequences of  points 
$\{ p_{1,n}\}$, $\{p_{2,n}\}$ in $\overline {\partial ^* \Omega}\cap H _0 ^+$  
  which  for every fixed $n$ are  distinct, with  the same projection onto $H_0$, and 
at infinitesimal distance from $H _ 0$ as $n \to + \infty$. 
\end{lemma}

\proof We argue by contradiction.
 Setting $t_{i, n} : =  {\rm dist} ( p  _{i,n} , H_0)$, we can assume up to a subsequence  that
$t _ {1, n} > t _{ 2 ,n }$ for every $n$.
We are going to show that
\begin{equation}\label{f:liminf} \liminf _{n \to + \infty} \frac{ \int_ \Omega h _{ p_{1,n}}    -\int_ \Omega h _{ p_{2,n}}    }
 {t _{ 1 ,n}- t _{ 2 ,n} } >0 \ ,
\end{equation} 
against the fact that $\Omega$ is $h$-critical. 

By the equality \eqref{f:cake1} in Lemma \ref{l:cake}, we have
\begin{equation}\label{f:minus} \begin{array}{ll} \displaystyle & \displaystyle \int_ \Omega h _{ p_{1,n}}    -\int_ \Omega h _{ p_{2,n}}   
 \displaystyle  = \int _{\sigma}  ^ { + \infty} \!\!\!  |\Omega \cap {B _ { r (s) } (p_{1,n} )} | - 
|\Omega \cap {B _ { r (s)}  (p_{2,n} )} |
   \, ds \\   
   \noalign{\medskip}
&   \displaystyle =  \int _{\sigma}   ^ {+ \infty}\!\!\!  |\Omega \cap (B _ {r(s)} (p_{1,n})  \setminus B _{ r(s)} (p_{2,n}) ) | - | \Omega \cap (B _ {r(s)}  (p_{2,n})  \setminus B _ {r(s)} (p_{1,n})) | \, ds \, . \end{array}
\end{equation} 
Since $\Omega$ is  not $h$-degenerate, and recalling the equality \eqref{f:cake2} in Lemma \ref{l:cake}, there exists a positive constant $C$ such that
\begin{equation}\label{f:plus} \frac{ \int _\Omega  | h _ { p_{1,n}} - h _ { p_{2,n}}| \, dy }  {t _{ 1 ,n}- t _{ 2 ,n} } = \frac{\int _{\sigma}   ^ {+ \infty}  
\big | \Omega  \cap
\big ( B _ { r (s) } (p_{1,n}) \Delta B _ { r (s) } (p_{2,n}) \big )  
 \big   | \, ds } {t _{ 1 ,n}- t _{ 2 ,n} }  \geq C \,.
\end{equation}
In view of \eqref{f:minus} and \eqref{f:plus},  the inequality \eqref{f:liminf} holds true provided, for $n$ large enough,
\begin{equation}\label{f:goal}
\frac{\int _{\sigma}   ^ {+ \infty} 
| \Omega \cap (B _ {r(s)}  (p_{2,n})  \setminus B _ {r(s)} (p_{1,n})) | \, ds} {t _{ 1 ,n}- t _{ 2 ,n} }
 \leq \frac{C}{4} \,.
\end{equation} 
By the inclusion $\Omega \subset H _0 ^+$, to prove the inequality \eqref{f:goal} it is sufficient to have \begin{equation}\label{f:goal2}
{ \int _{\sigma}   ^ { + \infty}  | H_0 ^+\cap (B _{ r(s)} (p_{2,n})  
 \setminus B _ {r(s)} (p_{1,n})) | \, ds }   = o ( {t _{ 1 ,n} - t _{ 2 ,n}} )    \,.\end{equation}

It remains to prove \eqref{f:goal2}. To that aim,
we  are going to provide two distinct estimates valid for $n$ large enough for the integrand in \eqref{f:goal2}, according to the values of the radius $r ( s)$.  More precisely, we distinguish the two  regimes  
$$
r (s) \leq k _n t _{1, n} \quad \text{ and } \quad  r ( s) > k _n t _{1, n} \,,$$
where $k_n$  is a constant larger than $1$, which will be suitably chosen at the end of the proof. 
In the estimates below, we  set $\gamma _n :=   {t _{ 1 ,n} - t _{ 2 ,n}}$; moreover, 
 we omit for shortness the index $n$, by
simply writing $p_1$, $p _2$, $t_1$, $t_2$, and $\gamma$.

\smallskip

\begin{itemize}
 \item{}

For $r (s) \leq  k t _ 1$,  we  have
 $$ | H _0 ^ + \cap ( B _{ r(s)} (p_{2})  \setminus B _ {r(s)} (p_{1}))| \leq |B _{ r(s)} (p_{2})  \setminus B _ {r(s)} (p_{1})| \leq  \omega _{d-1}   r ( s) ^ { d-1} \gamma  \,.$$

 \smallskip
\item{} For $r (s) >  k t _ 1$, 
 since $k > 1$ both the balls $B _ {r(s)} (p_{1})$ and $B _{ r(s)} (p_{2})$ intersect $H _ 0$.  
 Let us denote by $z (= z_n)$ the common projection of $p  _{1}$ and $p  _{2}$  onto $H_0$. 
The measure of $H _ 0 ^ + \cap (B _{ r(s)} (p_{2})  \setminus B _ {r(s)} (p_{1}))$
is not larger than the measure of the  region $D (s)$ obtained as the difference between two right cylinders having as axis the perpendicular to $H _ 0$ through $z$, as bases the $(d-1)$-dimensional ball contained into $H_0$ 
centred at $z$ with
radii  $r_2:= ( r(s) ^ 2 - t_{2} ^ 2  ) ^ {1/2}$ and  $r_1:= ( r(s)^ 2 - t_{1} ^ 2 ) ^ {1/2}$, and as heigh  $t _ 2 + \frac{ \gamma}{2}$. We have
$$|D (s) | = \omega _ { d-1} ( r_2 ^ {d-1} - r _ 1 ^ { d-1})\,  \big (t _ 2 + \frac{\gamma}{2} \big ) 
   \,.$$ 
By the convexity of the map $t \mapsto t ^ { d-1}$ (for $d \geq 2$) we infer
$$
\begin{array}{ll}  |D (s) | & \displaystyle  \leq (d-1) \omega _{d-1}   r _ 2 ^ { d-2} ( r _ 2 - r _ 1) \,  \big (t _ 2 + \frac{\gamma}{2} \big ) 
\\  \noalign{\medskip} 
\displaystyle  & = (d-1) \omega _{d-1}   r _ 2 ^ { d-2} (  \sqrt{ r(s) ^ 2 - t_{2} ^ 2  }   - \sqrt { r(s)^ 2 - t_{1} ^ 2 }  ) \, \,  \big (t _ 2 + \frac{\gamma}{2} \big ) 

\\  \noalign{\medskip} 
 & \displaystyle  = (d-1) \omega _{d-1}   r _ 2 ^ { d-2} \frac {t_1 ^ 2 - t _ 2 ^2 } {  r_1       + r_2  }   \,  \,  \big (t _ 2 + \frac{\gamma}{2} \big ) 
\\  \noalign{\medskip} 
 & \displaystyle  \leq  2  (d-1) \omega _{d-1}    r _ 2 ^ { d-3 } t _ 1  \,   \big (t _ 2 + \frac{\gamma}{2} \big )  \, \gamma \,,
\\  \noalign{\medskip} 
 & \displaystyle  \leq {  2  (d-1) \omega _{d-1}     }   r  (s)^ { d-3}  \, t _ 1 ^ 2    \,  \gamma 
 \\  \noalign{\medskip} 
 & \displaystyle  \leq {   \frac{  2  (d-1) \omega _{d-1}     }{k ^ 2} }   r  (s)^ { d-1}      \,  \gamma 
   \,. $$ 

\end{array}
 $$
\end{itemize} 
Now, if we set
$$s \big  (\lambda \big )  : = \sup \{ s :  r ( s)  >  \lambda  \}\,, \qquad \forall \lambda >0 \,,
$$ 
for $n$ sufficiently large  it holds 
$$\sigma = s ( { \rm diam} \Omega) <  s (k  t_1  )   \,.$$ 
Hence, 
$$\begin{array}{ll} 
& \displaystyle  \int _{\sigma}   ^ { + \infty}  | H_0 ^+\cap (B _{ r(s)} (p_{2})  \setminus B _ {r(s)} (p_{1})) | \, ds  = 
 \\ \noalign{\medskip} 
&  \displaystyle  \int _{s (k  t_1  ) }   ^ { + \infty }  | H_0 ^+\cap (B _{ r(s)} (p_{2})  \setminus B _ {r(s)} (p_{1})) | \, ds    +      \int _{\sigma }   ^ {  s (k t_1  )}  | H_0 ^+\cap (B _{ r(s)} (p_{2})  \setminus B _ {r(s)} (p_{1})) | \, ds     
   \leq  \\ \noalign{\medskip} 
& \displaystyle         \Big [ 
\omega _{ d-1}  \int _{s (k  t_1  )  }   ^ {+ \infty }   r(s) ^ {d-1}  \, ds     +  
\frac{  2  (d-1) \omega _{d-1}     }{k ^ 2}  \int _{\sigma} ^ {s (k  t_1  )  }     r  (s) ^ { d-1 }        \, ds         \Big ] \, \gamma 

\,.
     \end{array} $$  
     Finally  we claim that, by choosing $k= k _n \to + \infty $ in such way that $k  t_1 \to 0$,  the two addenda in square bracket are infinitesimal.  To prove such claim it is enough to have
     \begin{equation}\label{f:finito} 
     \int _{\sigma} ^ { + \infty}  r(s) ^ {d-1}  \, ds    < + \infty\,.
     \end{equation}
     Indeed in this case the first addendum will be  infinitesimal because $s ( k t _ 1)$ tends to $+ \infty$, while the second one
     will be infinitesimal because it is  bounded from above by a finite integral times a ratio which tends to zero. 
Eventually, condition \eqref{f:finito} holds true since  the convergence of the integral near $+ \infty $ is guaranteed by assumption \eqref{f:integrability},  while 
the convergence near $\sigma$ is guaranteed by the fact that, if $\sigma = 0$,  we have $r ( 0) \leq {\rm diam} \, \Omega$.

\qed

\bigskip 
Assume now that Claim 1 is false. Then, 
 either there exists 
$\{t _n\}\to 0$ such that $\forall n$  $\Omega _{t_n}  \cup  \mathcal R _{t_n} $ is {\it not} Steiner symmetric about $H _ { t_n}$, 
or there exists $\{t _n\} \to 0$ such that $\forall n$  $|\mathcal R _{t_n}  \setminus \Omega| >0$. 
Lemma  \ref{l:nopairs} ensures that none of these two cases  is possible (the detailed contradiction argument works as in the proof of  Step 1 in \cite[Theorem 1]{BF}).

   \bigskip
\subsection{ Step 2 (away contact at the stopping time)} 
We set 
$$T:= \sup \Big \{ t >0 \ :\ \text{ for all  $s \in [0, t)$,  symmetric inclusion occurs without away contact} \Big \} \, . $$ 
Since $\Omega$ is bounded, we have $T < + \infty$. Then  we prove the following claims:
\medskip

$\bullet$ {\it Claim 2a.  Symmetric inclusion  holds at $T$
 with away contact or with close contact.}
{\medskip}

$\bullet$  {\it Claim 2b. Symmetric inclusion cannot hold with close contact and no away contact.} 

\bigskip
In order to prove these claims, \EEE
 we need to establish preliminarily that  symmetric inclusion without away contact  implies the
 ``away inclusion properties" stated in the next lemma (precisely, in \eqref{f:inclusion1} and \eqref{f:inclusion2}).  Below,  for any $\delta >0$ and $s\geq 0$, we set
$$
 U _{T- \delta} ^ s :=   \Big \{ x + ( 2 \delta + 2 s ) \nu \ :\ x \in \mathcal R _{T- \delta}  \Big \} \,.
$$
Moreover, 
we denote by 
  $E \oplus B _ R$ the collection of points of $\R ^d$ with distance less than $R$ from a set $E$.

\begin{lemma}  Assume that symmetric inclusion occurs without away contact at $T$, and let $\delta >0$ be fixed. Then: 
\begin{itemize}
\item[-- ]
 There exists $s _\delta >0$ such that, for every $s \in [0, s _\delta]$,   
 \begin{equation}\label{f:inclusion1} 
 U _{T- \delta} ^ s \subset \Omega\,.
 \end{equation}

\item[-- ] 
There exists $\eta=\eta_\delta$ such that 
\begin{equation}\label{f:inclusion2}
(U_{T-\delta}^0\oplus B_\eta) \cap H _{T+ \delta  } ^+  \sq \Om\,.
\end{equation} 
\end{itemize} 

\end{lemma}

\proof 
We argue by contradiction. 
 If the inclusion \eqref{f:inclusion1} was false, we could find an infinitesimal sequence $\{s_n\}$ of positive numbers, and a sequence  of points $\{x'_n\}$  of density $1$ for $U _{T- \delta} ^ {s_n}$  but of density $0$ for $\Omega$. Up to a subsequence, there exists $x' := \lim _n x' _n$. 
 By construction, we have  
 $x' \in \{ x + 2 \delta \nu : x \in \overline { \partial ^ *  \mathcal R _{T- \delta}} \}  \subset  \overline { \partial ^ * \mathcal R _T }$. 
But, since we are assuming that symmetric inclusion occurs without away contact at $T$, it is readily checked that
 $\overline { \partial ^ *  \mathcal R _T} \subseteq {\rm int} (\Omega ^ { (1)} )$. 
 Then $x' \in  {\rm int} (\Omega ^ { (1)} )$, against the fact that $x' _n$ are points of density $0$ for $\Omega$. 
 
In a similar way, 
if the inclusion \eqref{f:inclusion2} was false,  we could find a sequence $\{x_n\} \in (\Om^c)^{(1)}$  such that $ x_n \in  (U_{T-\delta}^0\oplus B_{\frac 1n}) \cap H _{T+ \delta} ^+$. Since $U_{T-\delta}^0 \sq \Om$ is open (as a consequence of Proposition 13 in \cite{BF}\EEE ), we could also find $y_n \in  (U_{T-\delta}^0\oplus B_{\frac 3n}) \cap H _{T+ \delta_0 } ^+$ such that $y_n \in \partial ^* \Om$ (otherwise by Federer's Theorem, the perimeter of $\Om$ inside the set $B_{\frac 2n} (y_n)$ would be zero). By compactness, we would obtain  a limit point of  the sequence $\{y_n\}$ lying both in $\overline{ \partial ^* \Om}$ and in $\overline U_{T-\delta}^0 \cap \overline H _{T+ \delta_0  } ^+$, in contradiction with our assumption of symmetric inclusion without away contact at $T$. 
\qed

\medskip
\bigskip

{\bf Proof of Claim 2a}. The same arguments used to obtain the homonym claim in the proof of \cite[Theorem 1]{BF} apply.  

 \medskip 
 \bigskip
 {\bf Proof of Claim 2b.} 
 
\begin{figure} [h] 
\centering   
\def\svgwidth{7cm}   
\begingroup%
  \makeatletter%
  \providecommand\color[2][]{%
    \errmessage{(Inkscape) Color is used for the text in Inkscape, but the package 'color.sty' is not loaded}%
    \renewcommand\color[2][]{}%
  }%
  \providecommand\transparent[1]{%
    \errmessage{(Inkscape) Transparency is used (non-zero) for the text in Inkscape, but the package 'transparent.sty' is not loaded}%
    \renewcommand\transparent[1]{}%
  }%
  \providecommand\rotatebox[2]{#2}%
  \ifx\svgwidth\undefined%
    \setlength{\unitlength}{265.53bp}%
    \ifx\svgscale\undefined%
      \relax%
    \else%
      \setlength{\unitlength}{\unitlength * \real{\svgscale}}%
    \fi%
  \else%
    \setlength{\unitlength}{\svgwidth}%
  \fi%
  \global\let\svgwidth\undefined%
  \global\let\svgscale\undefined%
  \makeatother%
  \begin{picture}(3, 0.8)%
    \put(0.5 , 0 ){ \includegraphics[height=6cm]{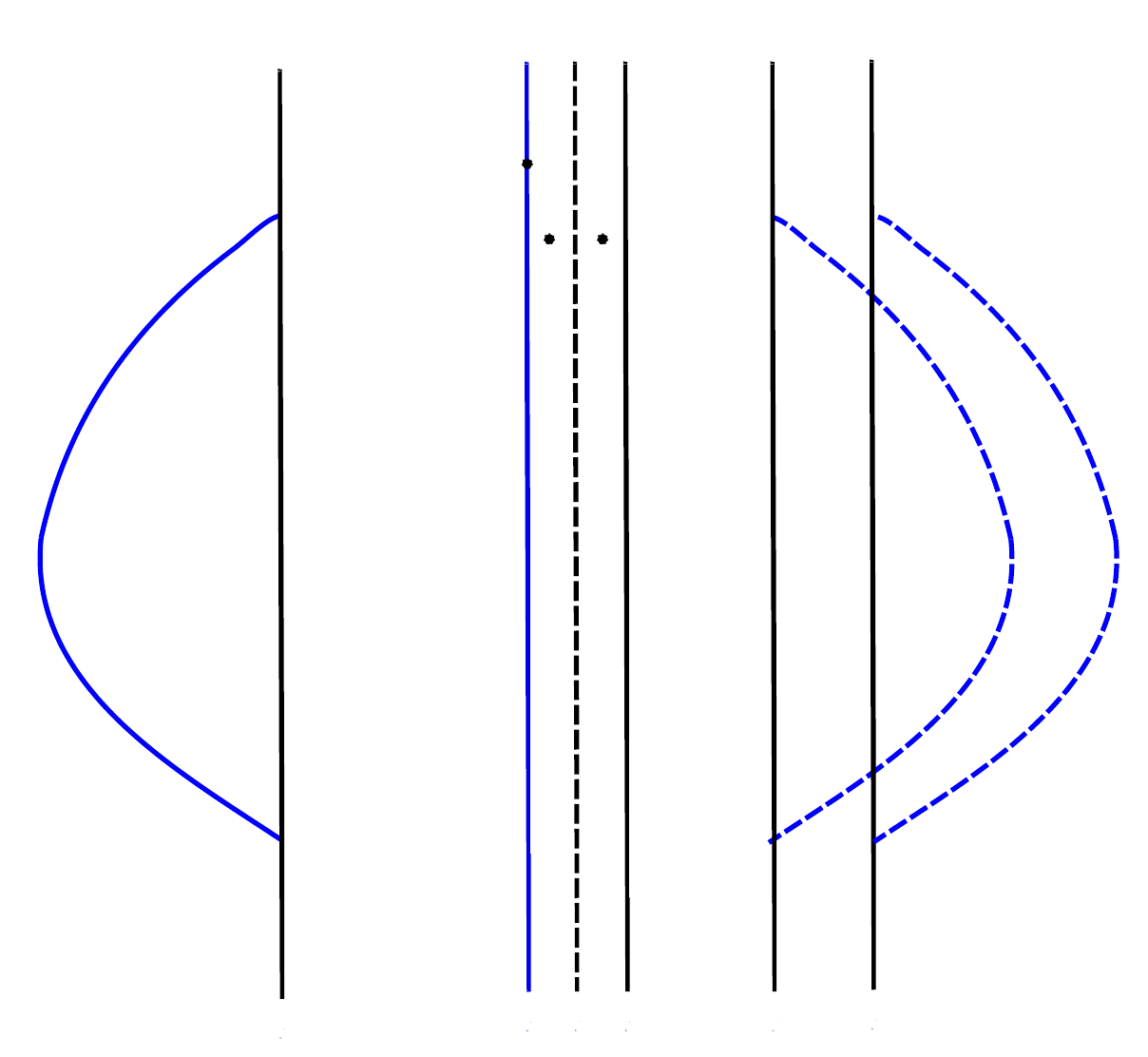}} %
    \put(0.587, 0.4){\color[rgb]{0,0,0}\makebox(0,0)[lb]{\small \smash{$\Omega _{T - \delta}$}}}
    \put(1.225, 0.4){\color[rgb]{0,0,0}\makebox(0,0)[lb]{\smash{\small $U^0 _{T - \delta}$}}}
    \put(1.343, 0.4){\color[rgb]{0,0,0}\makebox(0,0)[lb]{\smash{\small $U^{\tau _n} _{T - \delta}$}}}
 \put(0.67,  0.82){\color[rgb]{0,0,0}\makebox(0,0)[lb]{\smash{\small {$T- \delta$}}}}    
 \put(0.93,  0.82){\color[rgb]{0,0,0}\makebox(0,0)[lb]{\smash{\small {$T$}}}}     
\put(1.08,  0.82){\color[rgb]{0,0,0}\makebox(0,0)[lb]{\smash{\small {$T+ \delta$}}}}    
        \put(0.98,  -0.01){\color[rgb]{0,0,0}\makebox(0,0)[lb]{\smash{\small {$T\!\! +\!\! 2 \tau_n $}}}} 
           \put(1.17,  -0.01){\color[rgb]{0,0,0}\makebox(0,0)[lb]{\smash{\small {$T\!\! +\!\! \delta+\!\! 2 \tau_n  $}}}} 
    \end{picture}%
\endgroup%
\caption{The geometry in the proof of Claim 2b: the vertical lines represent, starting from the left, the  hyperplanes 
$H _ {T- \delta}$, $H _ T$ (in blue), $H_{T + 2 \tau _n}$, $H _ {T + \delta}$, and $H _ {T _ {\delta + 2 \tau _n}}$. }\label{figura}   
\end{figure}
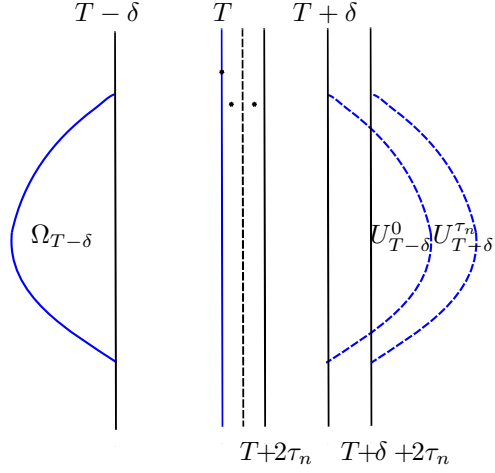

Assume that symmetric inclusion holds at $T$ with close contact. We are going to  contradict \eqref{f:hyp3} by showing that, if
$\{q_{1,n}\}$ and
  $\{q _{2,n}\}$ are sequences converging to a point $q\in H _ T$ as in \eqref{f:touching2}, 
it holds 
\begin{equation}\label{f:liminf2} \liminf _{n \to + \infty} \frac{\int_\Omega h _{q_{1,n}}  - \int_\Omega h _{q_{2,n}} }
{ \|  q _{1, n}-q _ { 2, n} \|}  >0 \,.
\end{equation} 

Below, we set for brevity 
$$\gamma _n :=    \|  q _{1, n}-q _ { 2, n} \| \,, $$ 

We fix $\delta >0$ small enough (to be chosen later), and  we let 
$\tau _n >0$ be such that $H _ { T + \tau _n }$ contains the midpoint of the segment $(q_{1,n}, q_{2,n})$ (see Figure \ref{figura}).  Up to working with $n$ sufficiently large, since $q_{1,n}$ and  $q_{2,n}$ converge to a point of $H _ T$, thanks to the away inclusion property \eqref{f:inclusion1}  we can assume that  \begin{equation}\label{f:contenimento} 
 U _{T- \delta} ^ {\tau _n}  \subset \Omega \,.
 \end{equation}
 We can also assume that  
 \begin{equation}\label{f:minoration}
 \gamma_n+ {\rm dist} \, (q_{1,n}, H_T)< \frac \delta 4\,.
 \end{equation}

We decompose
$$ \int_ \Omega h _{ q_{1,n}}    - h _{ q_{2,n}}   = X_n(\delta)+Y_n(\delta)+Z_n(\delta)$$
where
$$X_n(\delta)= \int_{ \Omega \cap  (     H _{T+ \tau_n} \oplus B _ {\delta + \tau _n}  ) } ( h _{ q_{1,n}}    - h _{ q_{2,n}})  $$
$$Y_n(\delta)=\int_{ \Omega_{T-\delta} \cup U_{T-\delta}^{\tau_n}   } ( h _{ q_{1,n}}    - h _{ q_{2,n}} )  $$
$$Z_n(\delta)=\int_{ \Omega \cap     H _{T+ \delta+ 2\tau_n} ^+\setminus U_{T-\delta}^{\tau_n} } ( h _{ q_{1,n}}    - h _{ q_{2,n}}) .$$
By symmetry, it holds
 $$Y_n(\delta)=0 \,.$$
From the non-degeneracy hypothesis we know that
\begin{equation}\label{f:nondegC} C:= \inf _n \frac{ \int_ \Omega |h _{ q_{1,n}}    - h _{ q_{2,n}} | }{\gamma_n} >0  \,. 
\end{equation}
In order to prove  \eqref{f:liminf2} we are going to  show first that
\begin{equation}\label{f:stimaX}
|X_n (\delta) |   \leq 
C  (\delta)   \gamma_n  \,, \quad \text{ with } C ( \delta) = o ( 1) \text{ as } \delta \to 0\,
\end{equation}
so that in particular we can choose $\delta$ so small that $C ( \delta) < \frac{C}{16}$, for $C$ as in \eqref{f:nondegC},  
and second that  
\begin{equation}\label{f:stimaZ}
Z_n(\delta) \ge C' \gamma_n \quad \text{ for some }  C'> \frac {C}{16} \,. 
\end{equation}

{\bf Proof of \eqref{f:stimaX}}. 
We decompose 
$$ \int_ \Omega |h _{ q_{1,n}}    - h _{ q_{2,n}} |  = X'_n(\delta) +Y'_n(\delta)+Z'_n(\delta),$$
where
$$X'_n(\delta)= \int_{ \Omega \cap  (     H _{T+ \tau_n} \oplus B _ {\delta + \tau _n}  ) } |h _{ q_{1,n}}    - h _{ q_{2,n}} | $$
$$Y'_n(\delta)=\int_{ \Omega_{T-\delta} \cup U_{T-\delta}^{\tau_n}   } |h _{ q_{1,n}}    - h _{ q_{2,n}} | $$
$$Z'_n(\delta)=\int_{ \Omega \cap     H _{T+ \delta+ 2\tau_n} ^+\setminus U_{T-\delta}^{\tau_n} } |h _{ q_{1,n}}    - h _{ q_{2,n}} |.$$

Assume by a moment to know that there exists a positive constant $C= C (\delta)$ such that 
\begin{equation}\label{f:eff1} X'_n (\delta)    \leq 
C  (\delta)   \gamma_n  \,, \quad \text{ with } C ( \delta) = o ( 1) \text{ as } \delta \to 0\,. 
\end{equation}

If this is the case,  since $|X_n|\le X'_n$, \eqref{f:stimaX} is satisfied. In particular, we can assume that $\delta$ is fixed such that $X'_n(3\delta) \le \frac {C}{16} \gamma_n$ for $n$ large enough.
Note that $X'_n$ is non decreasing in $\delta$,  so that $X'_n(\delta) \le \frac {C}{16} \gamma_n$.  

Since $|X_n|\le X'_n$, we have $|X_n(3\delta)| \le \frac {C}{16} \gamma_n$ and $|X_n(\delta)| \le \frac {C}{16} \gamma_n$.

We now prove \eqref{f:eff1}. To that aim, we are going to proceed in a similar way as done in the proof of Lemma \ref{l:nopairs}.  We set
$$\begin{array}{ll} 
& \displaystyle t _{ 1, n}:=  {\rm dist} ( q _ { 1, n}, H _ { T - \delta})  = \delta + \tau _n - \frac{\gamma _ n}{2}  
\\ \noalign{\medskip}  
&  \displaystyle t _ { 2, n} :=   {\rm dist} ( q _ { 2, n}, H _ { T - \delta})  = \delta + \tau _n + \frac{\gamma _ n}{2}  \,,
\end{array}$$

We have
$$X' _n (\delta) = 
\int _{\sigma} ^ {+ \infty}\!\!\!    | \Omega \cap  (     H _{T+ \tau_n} \oplus B _ {\delta + \tau _n}  ) \cap (B _ {r(s)} (q_{2,n})  \setminus B _ {r(s)} (q_{1,n}) ) |  \, ds\,.$$

Let us  provide two distinct estimates valid for $n$ large enough for the above integrand in the two   regimes
$$ r (s) \leq  \gamma _n + \sqrt \delta \ \quad \text{ and } \quad \ r ( s) > \gamma _n + \sqrt \delta   \,.$$

In the estimates below, we omit for shortness the index $n$, by
simply writing $q_1$, $q _2$, $t_1$, $t_2$, $\tau$, and $\gamma$.

\medskip

\begin{itemize}

\item{} For $r (s) \leq  \gamma + \sqrt \delta $,  we  have: 
$$
\begin{array}{ll}\displaystyle  |\Omega \cap  (     H _{T+ \tau} \oplus B _ {\delta + \tau }  ) \cap (B _ {r(s)} (q_{2})  \setminus B _ {r(s)} (q_{1}) )| &\displaystyle  \leq 
| B _{ r(s)} (q_{2})  \setminus B _ {r(s)} (q_{1}) | 
\\ \noalign{\medskip} & \leq 
\omega _ { d-1}  r (s) ^ { d-1}   \gamma \ .
\end{array}
$$

 \smallskip
\item{} For $r (s) >  \gamma + \sqrt \delta$, 
since for $n$ large enough  $t _ 1 < \gamma + \sqrt \delta  $, both the balls $B _ {r(s)} (q_{1})$ and $B _{ r(s)} (q_{2})$ intersect $H _ {T- \delta} ^+ $.  
 We have 
  $$ \begin{array}{ll} &  \displaystyle |   \Omega \cap  (     H _{T+ \tau} \oplus B _ {\delta + \tau }  ) \cap (B _ {r(s)} (q_{2})  \setminus B _ {r(s)} (q_{1}) ) |  \leq \\ 
  \noalign{\medskip} 
  & \displaystyle |  H _ {T- \delta} ^+  \cap (B _ {r(s)} (q_{2})  \setminus B _ {r(s)} (q_{1}) ) |  \,.
  \end{array}
   $$  
Let us denote by $z ( = z_n)$ the common projection of $q  _{1}$ and $q  _{2}$  onto $H_{T- \delta}$.  
The measure of $H _ 0 ^ + \cap (B _{ r(s)} (q_{2})  \setminus B _ {r(s)} (q_{1}))$
is not larger than the measure of the  region $D (s)$ obtained as the difference between two right cylinders having as axis the perpendicular to $H _ 0$ through $z$, as bases the $(d-1)$-dimensional ball contained into $H_0$ with centre $z$ and radii  $r_2:= ( r(s) ^ 2 - t_{2} ^ 2  ) ^ {1/2}$ and  $r_1:= ( r(s)^ 2 - t_{1} ^ 2 ) ^ {1/2}$, and as heigh  $t _ 2 + \frac{ \gamma}{2}$. As in the proof of Lemma \ref{l:nopairs}, we have
$$
|D (s) | 
 \leq    2  (d-1) \omega _{d-1}    r (s) ^ { d-3 } t _ 1 ^2 \, \gamma  \,,
 $$ 
Since, for $n$ large enough, we have $t _1 \leq  {\delta} \leq r ( s) \sqrt \delta$, we infer that

$$  | D (s) |    \leq 2  (d-1) \omega _{d-1}    r  (s) ^ { d-1 } \delta  \, \gamma  \,.   $$

\end{itemize}

If $s \big  (\lambda \big ):=  \sup \{ s :  r ( s)  >  \lambda  \}$, 
for large  $n$  it holds 
$\sigma = s ( { \rm diam} \Omega) <  s (\gamma + \sqrt \delta  )$, and hence 
$$\begin{array}{ll} 
& \displaystyle \int _{\sigma}   ^ { + \infty}  | H_0 ^+\cap (B _{ r(s)} (q_{2})  \setminus B _ {r(s)} (q_{1})) | \, ds     =
 \\ \noalign{\medskip} 
& \displaystyle \int _{s (\gamma + \sqrt \delta  ) }   ^ { +\infty }  | H_0 ^+\cap (B _{ r(s)} (q_{2})  \setminus B _ {r(s)} (q_{1})) | \, ds    + 
\int _{\sigma }   ^ {  s (\gamma + \sqrt \delta  )}  | H_0 ^+\cap (B _{ r(s)} (q_{2})  \setminus B _ {r(s)} (q_{1})) | \, ds   \leq  
\\ \noalign{\medskip} 
& \displaystyle        \Big [  
\omega _ { d-1} \int _{s (\gamma + \sqrt \delta  )  }   ^ { + \infty }   r(s) ^ {d-1}  \, ds    +  2  (d-1) \omega _{d-1}  \delta  \int _{\sigma} ^ {s (\gamma + \sqrt \delta  )  }     r  (s) ^ { d-1 }        \, ds 
    \Big ]  \, \gamma 
\,.
     \end{array} $$  
Then the proof of \eqref{f:eff1} is achieved  by taking $C ( \delta)$ equal to 
the expression in square bracket. Indeed,  recalling \eqref{f:finito}, we see that such expression is infinitesimal as $\delta \to 0$:   the  first addendum
     is  infinitesimal because $s ( \gamma + \sqrt \delta)$ tends to $+ \infty$, while the second one
     is  bounded from above by a multiple of $\delta$.

\EEE

\bigskip
{\bf Proof of \eqref{f:stimaZ}}. 
Since on $H _{T+ \delta+ 2\tau_n} ^+$ the difference $ h _{ q_{1,n}}    - h _{ q_{2,n}}$ is non negative, we get that $Z_n(\delta)=Z'_n(\delta)$. There are two situations:
$$\mbox{either $Y'_n(\delta)\ge Z'_n(\delta)$,\ \ \   or $Z'_n(\delta) > Y'_n(\delta)$}.$$
 If $Z'_n (\delta)>Y'_n(\delta)$, then from the non degeneracy hypothesis together with the estimate of $X'_n(\delta)$ (see respectively \eqref{f:nondegC} and \eqref{f:eff1}),  we get $$Z'_n(\delta) \ge \frac C4 \gamma_n$$
hence we get \eqref{f:stimaZ}. 

Assume now that $Y'_n(\delta)>Z'_n(\delta)$.  By arguing as above, we see that $Y'_n(\delta) \ge \frac C4 \gamma_n$. We shall prove that there exists a constant $K(\delta) >0$ such that  $Z'_n (\delta) \ge K(\delta)  \gamma_n$. Although we are not able to evaluate $K(\delta)$, its strict positivity will be  sufficient. In fact, following from our geometric construction and in particular from the definition of $Z' _ n (\delta)$, the map $\delta \mapsto K(\delta)$ is non increasing. Consequently,   if one replaces $\delta$ with a smaller value $\tilde \delta  <\delta$, it holds $K(\tilde \delta)\ge K(\delta)$ so that $Z'_n(\tilde \delta)   \ge K(\tilde \delta)  \gamma_n\ge K(\delta) \gamma_n$. In order to conclude the proof, it is then enough to fix $\tilde \delta $ small enough such that $X_n(3 \tilde \delta) \le \frac {K(\delta)}{4} \gamma_n$ and reproduce the same reasoning with $\tilde \delta$ instead of $\delta$. 

It remains to show that there exists a constant $K(\delta)>0$ such that  $Z'_n (\delta) \ge K(\delta) \gamma_n$. We know that $Y'_n(\delta) \ge \frac C4 \gamma_n$ or, equivalently,
$$\int _{\sigma} ^   { + \infty} \!\!\!  | (\Omega_{T-\delta} \cup U_{T-\delta}^{\tau_n}  )  \cap (B _ { r(s)}  (q_{1,n})  \Delta B _ {r(s)} (q_{2 ,n})) |   \, ds  \ge \frac C4 \gamma_n.$$
By  symmetry,
$$\int _{\sigma} ^   { + \infty} \!\!\!  |   U_{T-\delta}^{\tau_n}   \cap (B _ { r(s)}  (q_{1,n})  \setminus B _ {r(s)} (q_{2 ,n})) |   \, ds  \ge \frac C8 \gamma_n.$$
Since $X'_n(3 \delta) \le  \frac {C}{16} \gamma_n$, then
\begin{equation}\label{bf.01}
\int _{ [{\sigma}, { + \infty}] }  \!\!\!  1_{r(s) \in [2\delta, diam \Om]} |   U_{T-\delta}^{\tau_n}   \cap (B _ { r(s)}  (q_{1,n})  \setminus B _ {r(s)} (q_{2 ,n})) |   \, ds >0.
\end{equation}

By \eqref{bf.01},  setting
$$A(n, \delta, s)=   U_{T-\delta}^{\tau_n}   \cap (B _ { r(s)}  (q_{1,n})  \setminus B _ {r(s)} (q_{2 ,n})),$$ 
the one dimensional measure of the set 
$$\mathcal S:= \Big \{ s \ :\   r(s) \in [2\delta, {\rm diam} \, \Om] \ \text{ and } \  |A(n, \delta, s)| >0 \Big\}$$ 
is strictly positive. 

We have
$$Z' _n (\delta) \geq  \int_{\sigma}^{+ \infty} 1_{|A(n, \delta, s)| >0}  |\Omega \cap    (  H _{T+ \delta+ 2\tau_n} ^+\setminus U_{T-\delta}^{\tau_n} ) \cap ( B _{ r (s) }  (q_{1,n} ) \setminus  B _{ r (s) }  (q_{2,n} ) ) |  \, ds \,.$$ 

Now, for any $s \in \mathcal S$, we want to estimate from below the measure of the set 
$$  \Omega \cap    (  H _{T+ \delta+ 2\tau_n} ^+\setminus U_{T-\delta}^{\tau_n} ) \cap ( B _{ r (s) }  (q_{1,n} ) \setminus  B _{ r (s) }  (q_{2,n} ) ) \,.$$

To that end, on the straight line  through $q_{1,n}$ and $q_{2,n}$,  we look at the segment  
$$\big \{ P_{H_T} (q_{1,n}) + t \nu \  :\ t \in [\delta, \delta + {\rm diam} \Om] \big \}  \,, $$
where   $P_{H_T} (q_{1,n})$ is the orthogonal projection of $q_{1,n}$ onto $H _ T$.   
This segment can contain neither interior  points of $U_{T-\delta}^{\tau_n}$ nor  points of its  essential boundary, otherwise we would have an away contact point. Consequently,  its distance to the closure of the set $U_{T-\delta}^{\tau_n}$ is strictly positive, say $\eta ' >0$.

There are two possibilities: either all the set $(B _ { r(s)}  (q_{1,n})  \setminus B _ {r(s)} (q_{2 ,n})) \cap H _{T+ \delta+ 2\tau_n} ^+$ is contained in $\Om$, or not. 

In the first situation, the set 
$$(B _ { r(s)}  (q_{1,n})  \setminus B _ {r(s)} (q_{2 ,n})) \cap B_{\eta '} (P_{H_T} (q_{1,n})+(r(s)+\tau_n)\nu)$$
is contained in $\Om\cap  H _{T+ \delta+ 2\tau_n} ^+$ but does not intersect $U_{T-\delta}^{\tau_n}$. Moreover,  thanks to \eqref{f:minoration}, there exists a constant $C(\delta, d)$ depending on $\delta$ and the dimension of the space such that 
$$|(B _ { r(s)}  (q_{1,n})  \setminus B _ {r(s)} (q_{2 ,n})) \cap B_{\eta '} (P_{H_T} (q_{1,n})+(r(s)+\tau_n)\nu)| \ge \gamma_n C(\delta, d) (\eta')^{d-1}.$$

In the second situation, since $(B _ { r(s)}  (q_{1,n})  \setminus B _ {r(s)} (q_{2 ,n}))\cap H _{T+ \delta+ 2\tau_n} ^+$ contains already points from $\Om$ (precisely from $ U_{T-\delta}^{\tau_n}$), it also contain points from  $\partial ^* \Om$. Since any such boundary point lies at distance at least $\eta$ from $ U_{T-\delta}^{\tau_n}$ (recall \eqref{f:inclusion2}), we can find a point $x_n$ such that $ B_{\frac \eta2} (x_n) \cap ((B _ { r(s)}  (q_{1,n})  \setminus B _ {r(s)} (q_{2 ,n}))\cap H _{T+ \delta+ 2\tau_n} ^+ \sq  \Om \sm U_{T-\delta}^{\tau_n}$. As before,  thanks to \eqref{f:minoration}, there exists a constant $C(\delta, d, {\rm diam} \Om)$, now depending on the diameter as well, such that
$$|B_{\frac \eta2} (x_n) \cap ((B _ { r(s)}  (q_{1,n})  \setminus B _ {r(s)} (q_{2 ,n}))\cap H _{T+ \delta+ 2\tau_n} ^+| \ge 
\gamma_n C(\delta, d, {\rm diam} \,  \Om) \eta^{d-1}.$$
Finally,
$$Z'_n (\delta) \ge \int_{\sigma}^{+\infty } 1_{|A(n, \delta, s)| >0} \gamma_n (C(\delta, d)\wedge C(\delta, d, {\rm diam} \,  \Om) ) (\eta\wedge \eta')^{d-1}ds \ge K(\delta) \gamma_n,$$
for some positive constant $K(\delta)$. 

If $K(\delta)> \frac C8$, the proof is achieved. Otherwise, we choose $\tilde \delta < \delta$ such that $X_n(3 \tilde \delta ) \le \frac{K(\delta)}{4} \gamma_n$. Reproducing the same arguments and using the fact that $K(\tilde\delta)\ge K(\delta)$, we conclude the proof.

  \qed

 \bigskip
 
\subsection{Step 3  (decomposition of $\Omega$ into symmetric and non-symmetric part).}  
 We show that $\Omega$ can be decomposed as 
$$
\Omega= \Omega^s\sqcup \Omega ^{ns}\, , 
$$ 
where  $\Omega ^ s$ is an {\it open set}  representing the  Steiner symmetric part of $\Omega$, given by 
$$
\Omega ^ s := \bigcup \Big \{ (p, p') \ :\ p' \text{ is an away contact point, \ $p$ is its symmetric about $H _ T$} \Big \} \,, 
$$
$(p, p')$ being the open segment with endpoints $p$ and $p'$, 
and $\Omega ^ { ns}:= \Omega \setminus \Omega ^ {s}$ represents the non-symmetric part. 
More precisely, we prove the following two claims:

\smallskip
$\bullet$ {\it Claim 3a.} 
{\it  If $p'$ is an away contact point and 
 $p$ is  its symmetric about about $H _ T$, 
\begin{eqnarray}&  |(B_{r(s)} (p') \setminus B_{r(s)} (p) ) \cap  \Om ^{ns}| =0 \ \text{ for a.e. } s \in(\sigma, + \infty)  ; & \label{f:diff-balls} 
 \\ \noalign{\smallskip} 
 &\exists \e >0 \ :\ \big | B _ \e (p' ) \cap (\Omega \setminus \mathcal R _ T) \big |  = 0 \,,  \text{ and hence $\Omega ^ s$ is open.} & \label{f:claim}
 \end{eqnarray} 
 }

$\bullet$ {\it Claim 3b. Denoting by $ \Omega ^ s _ i$  
the open connected components of $\Omega ^ s$, 
it holds 
\begin{eqnarray}
& \text{ 
 $\overline{\partial ^ * \Omega^ s_i} \cap (H_ T^ \pm  \setminus H _ T)$  are  connected sets;} 
  & \label{f:cc} 
\\ \noalign{\medskip} 
& \overline {\partial ^* \Omega ^s} \cap \overline {\partial ^* \Omega ^{ns}} \subset H _ T\,.
 & \label{f:studiobordi}
\end{eqnarray} }
\begin{remark}\label{f:lontananza} We point out that  \eqref{f:diff-balls} implies that every connected component $\Om^ s _i$ of $\Om ^ s$ satisfies 
\begin{equation}\label{f:separation2} {\rm dist} (\Om ^ s _ i, \Om ^ {ns}) \geq \sup _{ s \in(\sigma, + \infty) } r(s) = r ( \sigma) \,.
\end{equation}  
\end{remark}

\bigskip 
{\it Proof of Claim 3a.} 
Starting from the assumption of $h$-criticality, we obtain
$$\begin{array}{ll}
0 & \displaystyle  = \int_\Omega  h _ {p'} -  \int _ { \Omega _T  \cup  \reflexT  }  h _ {p'} - \int_\Omega  h _ {p} +  \int _ { \Omega _T  \cup  \reflexT  }  h _ {p'} 
\\ \noalign{\medskip} 
& \displaystyle = \int_ { \Omega \setminus ( \Omega _T  \cup  \reflexT) }   h _ {p'} - \int_ { \Omega \setminus ( \Omega _T  \cup  \reflexT) }   h _ {p} 
\\ \noalign{\medskip} 
& \displaystyle = \int _0 ^ { + \infty} \int_ { \Omega \setminus ( \Omega _T  \cup  \reflexT) } \big [ \chi _ { B _ {r(s)} (p')}  (y)  -  \chi _ { B _ {r(s)} (p)} (y)  \big ]  \, dy \, ds 
\\ \noalign{\medskip} 
& \displaystyle = \int _0 ^ { + \infty} \int_ { \Omega \setminus ( \Omega _T  \cup  \reflexT) }  \chi _ { B _ {r(s)} (p')}  (y) \big [ 1 -  \chi _ { B _ {r(s)} (p)} (y)  \big ]  \, dy \, ds
\\ \noalign{\medskip} 
& \displaystyle = \int _0 ^ { + \infty} \big |  (\Omega \setminus ( \Omega _T  \cup  \reflexT)  )  \cap  ( B _ {r(s)} (p') \setminus  B _ {r(s)} (p) )  \big |   \, ds
\\ \noalign{\medskip} 
& \displaystyle = \int _0 ^ { + \infty} \big |  (\Omega \setminus  \reflexT  )  \cap  ( B _ {r(s)} (p') \setminus  B _ {r(s)} (p) )  \big |  \, ds
\\ \noalign{\medskip} 
& \displaystyle = \int _{\sigma} ^ { + \infty} \big |  (\Omega \setminus  \reflexT  )  \cap  ( B _ {r(s)} (p') \setminus  B _ {r(s)} (p) )  \big |  \, ds
 \,, 
\end{array}$$ 
which proves \eqref{f:diff-balls}.

\smallskip
Let us prove \eqref{f:claim}. We claim that there exists $s _0 \in (\sigma, + \infty) $ such that 
\begin{equation} \label{f:case3} 0 < \big | \big [ B _{r(s_0)}  (p' ) \setminus {B _ {r(s_0)}   (p)} \big ] \cap \Omega \big | < \big |  B _ {r(s_0)}   (p' ) \setminus {B _ {r(s_0)}   (p)} \big | \,.  
\end{equation}
Indeed, 
 since $\Omega$ is not $h$-degenerate, by equality \eqref{f:cake2} in Lemma \ref{l:cake} we have
 $$\mathcal L ^ 1 \big (  \big \{ s> \sigma \ :\ | \Omega \cap (B_{r(s)} (p) \Delta B_{r(s)} (p') )  | > 0   \big \}  \big )   >0\,.$$

Hence, using also \eqref{f:diff-balls}, we can pick $s_0 \in (\sigma, + \infty)$ such that the left inequality in \eqref{f:case3} is satisfied and
\begin{equation}\label{f:scelta} 
|(B_{r(s_0)} (p') \setminus B_{r(s_0)} (p) ) \cap  \Om ^{ns}| =0\,. 
\end{equation}
We observe that, for such $s_0$, also the right
 inequality in \eqref{f:case3} is necessarily satisfied. Indeed,  if this is not the case, we have  $$ \big  | \big [ B _ {r(s_0)}  (p' ) \setminus {B _ {r(s_0)} (p)} \big ] \cap \Omega \big | = \big | B _ {r(s_0)} (p' ) \setminus  {B _ {r(s_0)} (p)}  \big |\,.
$$
In view of  \eqref{f:scelta}, this implies
that   $B _ {r(s_0)}  (p' ) \setminus  {B _ {r(s_0)} (p)} $ is contained into $\reflexT$, and hence $B _ {r(s_0)} (p ) \setminus  {B _ {r(s_0)}(p')} $ is contained into $\Omega _T$.  
Since  $\Omega _T  \cup \reflexT $  is Steiner-symmetric 
about  $H _T$, we obtain (via Fubini Theorem) that $p$ and $p'$  belong to ${\rm int} (\Omega ^ { (1)})$, contradicting the fact that  they belong to  $\overline {\partial ^* \Omega}$. 

\smallskip
Now we observe that  
\begin{equation}\label{f:brexit}  \exists y' \in \big [ B_{r(s_0)}  (p') \setminus \overline{B _{r(s_0)} } (p) \big ] \cap \overline{\partial ^* \Omega} \,. \end{equation}   
Indeed, if  \eqref{f:brexit} was false, $B_{r(s_0)}  (p') \setminus \overline{B _ {r(s_0)} }  (p)$ would be contained either into ${\rm int}  ( \Omega ^ { (1)} ) $ or into ${\rm int}(\Omega ^ { (0)})$, against \eqref{f:case3}. 
In view of \eqref{f:scelta},  the two sets $\Omega$  and $\reflexT$ have the same density at  every point of
$B_{r(s_0)}  (p') \setminus \overline{B _ {r(s_0)} }  (p)$, and hence 
$$\big [ B_{r(s_0)}  (p' ) \setminus \overline{ B _ {r(s_0)} }  (p) \big ]  \cap   \partial ^  * \Omega = \big [ B_{r(s_0)}  (p' ) \setminus \overline{ B _ {r(s_0)} } (p) \big ] \cap \partial ^* \reflexT \,;
$$
consequently, since the set $B_{r(s_0)}  (p') \setminus \overline{B _ {r(s_0)}  (p)}$ is open, we have 
\begin{equation}\label{f:same-closure-bdry} 
 \big [ B_{r(s_0)} (p') \setminus \overline{B _ {r(s)} }  (p) \big ] \cap \overline {\partial ^* \Omega}  = \big [ B_{r(s_0)}  (p') \setminus \overline{B _ {r(s_0)} } (p) \big ] \cap \overline {\partial ^* \reflexT}  \,.
\end{equation}

By \eqref{f:brexit} and \eqref{f:same-closure-bdry}, it turns out that the point $y'$ is itself an away contact point. 
Therefore, denoting by $y$ its symmetric about $H _ T$, 
in the same way as we obtained \eqref{f:diff-balls}, replacing the pair $p, p'$ by the pair $y, y'$, we obtain 
 \begin{equation}\label{f:diff-balls2} 
\big | \big [ B _ {r(\tau)}  (y' ) \setminus B _ {r(\tau)}  (y )  \big ] \cap \big ( \Omega \setminus \reflexT \big ) \big | = 0 \qquad  \text{for a.e.} \ \tau \in (\sigma, + \infty) \,. 
\end{equation} 
Since the set $B_{r(s_0)}  (p') \setminus \overline{B _ {r(s_0)} } (p)$ is open,  for every $\e>0$ sufficiently small the ball  $B _ \e ( y')$ is contained into $B_{r(s_0)} (p') \setminus \overline{B _{r(s_0)}  } (p)$, and hence
\begin{equation}\label{f:specchio}
\exists \e (s_0) >0 \ :\ B _{\e(s_0)} (p') \subset  \big [ B_{r(s_0)}  (y') \setminus \overline{B _ {r(s_0)} } (y) \big ]\,.
\end{equation} 
This achieves the proof of \eqref{f:claim} in case the equality in \eqref{f:diff-balls2} is satisfied at $\tau  =  s_0$. 
But, since \eqref{f:diff-balls2} holds merely for a.e. $\tau \in (0, + \infty)$, we have to refine the argument as follows.  
By \eqref{f:specchio}, for $\lambda$ sufficiently close to $r ( s_0)$, we have
\begin{equation}\label{f:specchio-largo}
p' \in   \big [ B_{\lambda}  (y') \setminus \overline{B _ {\lambda} } (y) \big ]\,.
\end{equation} 
Then, by the continuity from the right of the map $r \mapsto r (s)$, there exists $\delta >0$ such that
\begin{equation}\label{f:intornodx}
\forall s \in [s_0, s_0 + \delta)\,, \, \quad p' \in  \big [ B_{r(s)}  (y') \setminus \overline{B _ {r(s)} } (y) \big ]\end{equation} 
and hence
\begin{equation}\label{f:finale}
\forall s \in [s_0, s_0 + \delta)\,, \, \quad \exists \e (s) >0 \ :\ B _{\e(s)} (p') \subset  \big [ B_{r(s)}  (y') \setminus \overline{B _ {r(s)} } (y) \big ]\,.
\end{equation}
Eventually, the proof of  \eqref{f:claim} is achieved by choosing
$ s \in [s_0, s_0 + \delta)$ such that 
 the equality in \eqref{f:diff-balls2} is satisfied at $\tau  =  s$. \qed

\bigskip {\it Proof of Claim 3b.}   The same arguments used to obtain the homonym claim in the proof of \cite[Theorem 1]{BF} apply.   \EEE

 \subsection{Step 4 (conclusion)} 
We show that the open connected components of $\Omega ^ s$  are balls
of the same radius $R> \eta/2$, lying at distance larger than or equal to  $r (\sigma)$, while the set $\Omega ^ { ns}$ is Lebesgue negligible.  

 In order to   formulate more precisely the claims which conclude our proof, we need to set up some additional definitions and  notation.  \EEE
  
  \smallskip 
Given two different open connected components $\Om_i^{s} , \Om_j ^ {s}$ of   $\Om^s$,  we say that $\Om^s_i$ {\it is in   $h$-contact with $ \Om_j ^ s$} if   there exists an away contact point $p' \in \overline{\partial ^* \Om^s _i} \setminus H_T$  such that, denoting by $p$ its symmetric about $H _ T$, it holds
$$  
\int_{ \Om^s _{j}} \big | h _ p - h _ {p'}  \big |>0.
$$ 
It is not difficult to check that,   if $\Om^s_i$ is in   $h$-contact with $ \Om_j ^ s$,  $\Om^s_j$  is in   $h$-contact with $ \Om_i^ s$.  

If $\Omega _ i ^ s$ is not in contact with  any other component of $\Omega ^ s$, we say that $\Omega_ i ^ s$ is {\it $h$-isolated}.   

\begin{remark}\label{r:isolation} We observe that, if $\Omega^s _i$ is $h$-isolated, for every  $p,p' \in \overline{\partial ^* \Om^s _i} \setminus H_T$ symmetric about $H_T$, it holds 
$$  
\int_{ \Om^s } \big | h _ p - h _ {p'}  \big | =  \int_{ \Om^s _{i}} \big | h _ p - h _ {p'}  \big |  \,.
$$ 
Hence, if $\Omega ^ s _j$ is any other component of $\Omega ^ s$, we have 
$$\int _{\sigma} ^ { \infty} 
 |(B_{r(s)} (p') \setminus B_{r(s)} (p) ) \cap  \Om ^s _j | \, ds = 0 \, , $$  
 which implies 
\begin{equation}\label{f:separation}
{\rm dist} (\Omega ^ s _ j, \Omega ^ s _i) >   \sup _ { s \in(\sigma, + \infty) } r(s) = r ( \sigma) \,.
\end{equation} 
\end{remark} 

\bigskip
 Since our strategy will require to let the initial hyperplane vary, we will write
$$\Omega = \Omega ^{\nu, s} \sqcup \Omega ^ {\nu, ns} \,,
$$ where the additional superscript $\nu$ indicates the direction of the parallel movement, namely the normal to the initial hyperplane $H_0$ (and the decomposition is always meant with respect  to the parallel hyperplane $H _ T$ at the stopping time $T$ defined in Step 2). 
  
  \smallskip
The fourth and final step of our proof consists in showing the following claims:  \EEE 
    
  \smallskip 
  
  $\bullet$ 
{\it Claim 4a.  Given $\nu \in \mathbb S ^ { d-1}$,  let $\Omega _ \flat$ be a  $h$-isolated open connected component  of $\Omega ^ {\nu, s}$.
Then $\Omega _ \flat$ is a ball of radius at least $\eta/2$,  and  $\Omega \setminus \Omega _\flat$ is $h$-critical and  not $h$-degenerate, unless it has measure zero. }
      
   \smallskip
$\bullet$ {\it Claim 4b. The following family is empty: } 
  $$\mathcal F: = \bigcup _{\nu \in \mathbb S^ { d-1}}  \Big \{ \text{\it open connected components not $h$-isolated of }  \Omega ^ {\nu, s}  \Big \}\,. $$   

\smallskip
$\bullet$ {\it Claim 4c (conclusion).}  {\it $\Omega$ is equivalent to a finite union of balls of radius $R> \eta/2$, at mutual distance larger than or equal to $r ( \sigma) $.  } 
 
 \bigskip

 {\it Proof of claim 4a}.  Given $\nu \in \mathbb S ^ { d-1}$, let $\Omega _ \flat$ be a  $h$-isolated open connected component  of $\Omega ^ {\nu, s}$. 
Assume by a moment to know that 
\begin{equation}\label{f:grasse}
\Omega _ \flat \text{  is  $h$-critical and  not $h$-degenerate.}   
\end{equation}
In this case, we can restart our proof, with $\Omega _\flat$ in place of $\Omega$. 
Given an arbitrary direction $\widetilde \nu \in \mathbb S ^ { d-1}$,  we make the decomposition 
$$\Omega _\flat = \Omega _\flat  ^ {\tilde \nu, s} \sqcup \Omega _\flat  ^ {\tilde \nu, ns}\,.$$
  It is not difficult to show that, unless $\Omega _\flat ^ {\tilde \nu, ns}$ is empty,  this decomposition 
splits $\Omega _\flat$ into two open sets, contradicting the connectedness of $\Omega  _\flat$.

(The detailed argument can be found in the proof of \cite[Theorem 1]{BF}, see Claim 4a. therein).    \EEE
Hence   $\Omega _\flat$ is Steiner symmetric about a hyperplane with unit normal $\widetilde \nu$.  
By the arbitrariness of $\widetilde \nu$, we deduce that $\Omega _\flat$ is 
a ball. 
  Denote by $R$ the radius of this ball.   
  We observe that a ball  of radius $R$ is not $h$-degenerate if and only if $\mathcal L ^ 1 \big ( \big \{ s \, :\, r(s) \in (0, 2 R)  \} \big ) > 0 $. Recalling \eqref{f:misura0}, the nondegeneracy of $\Omega _\flat$ yields the lower bound $R >  \eta  /2$.

  \smallskip
  To conclude the proof of Claim 4a., it remains to show that \eqref{f:grasse} holds true and that the same property is valid for $\Omega \setminus \Omega _\flat$, unless it has measure zero.

  For the sake of clearness, we split the proof into 
  three consecutive lemmas.

    \begin{lemma}\label{bfm03} Under  the assumptions of Theorem \ref{t:serrin3}, 
given $\nu \in \mathbb S ^ { d-1}$, let $\Omega _ \flat$ be a  $h$-isolated open connected component  of $\Omega ^ {\nu, s}$. Then  
$$
 \inf_{x_1,x_2 \in \partial ^* \Om_\flat} \frac{\int_{\Om^{\nu, s}} | h_ {x_1}- h_{x_2}|}{\|x_1-x_2\|} >0.
 $$  \end{lemma}
  \proof 
   \bigskip
 Assume by contradiction that  \begin{equation}\label{f:hyp12}
 \inf_{x_1,x_2 \in \partial ^* \Om_\flat} \frac{\int_{\Om^{\nu, s}} | h_ {x_1}- h_{x_2}|}{\|x_1-x_2\|} =0.
\end{equation}  
Then there exist sequences of distinct points $\{x_{1,n}\}, \{x_{2,n}\}  \subset  \partial ^* \Om_\flat$, with
  $\|x_{1,n} -x_{2,n} \|\ra 0$,  such that  
 $$ \frac{\int_{\Om^{\nu, s}} | h_ {x_{1,n}}- h_{x_{2,n}}|}{\|x_{1,n}-x_{2,n}\|}\ra 0  \,.$$

   Up to subsequences, we may assume that $\|x_{1,n}-x_{2,n}\|$ converges to $0$ decreasingly,  and that $\{x_{1,n}\}$ and  $\{x_{2,n}\}$ converge to some point $\overline x \in \overline {\partial ^* \Om_\flat }$, which may belong or not to $H _ T$, being as usual $T$ the stopping time
    defined as in Step 2 for the parallel movement with normal $ \nu$. 
   Let us examine the two cases separately.

   In case $\overline x \not \in H _ T$, we may assume without loss of generality that  $\{x_{1,n}\}, \{x_{2,n}\} \subset H _ T ^ +\setminus H _ T$.

   Since $\Omega$ is not $h$-degenerate, \eqref{f:hyp12} implies 
$$
 \inf_{x_1,x_2 \in \partial ^* \Om_\flat} \frac{\int_{\Om^{\nu, ns}} | h_ {x_1}- h_{x_2}|}{\|x_1-x_2\|} >0.
$$

In particular, for $n=1$, we have   
$$ \int _{\sigma} ^ { + \infty} |\Om^{\nu, ns} \cap (B_{r(s)} (x_{1, 1}) \Delta B_{r(s)} (x_{2, 1})) |  \, ds =   \int_{\Om^{\nu, ns}} | h_ {x_{1,1}}- h_{x_{2,1} }|  >0\,.$$
 
We infer that there exists $s _0 \in (\sigma, + \infty)$  such that  
$$  |\Om^{\nu, ns} \cap (B_{r(s_0)} (x_{1, 1}) \Delta B_{r(s_0)} (x_{2, 1})) | >0 \,.$$ 
Hence 
 we can pick a point $p\in {\rm int}(B_{r(s_0)} (x_{1, 1}) \Delta B_{r(s_0)} (x_{2, 1}))$ of density $1$ for $\Omega ^ { \nu , ns} $, and   
   a radius $\e >0$  sufficiently small so that 
  \begin{equation}\label{f:alto} 
  | B_\varepsilon (p) \cap \Om^{\nu, ns} | \geq \frac{1}{2}  | B_\varepsilon (p)| \,.
  \end{equation}
Possibly reducing $\e$ we can also assume that 
 $B_\varepsilon (p)  \sq  \big (B_{r(s_0) } (x_{1, 1}) \Delta B_{r(s_0) } (x_{2, 1})\big )$.

  \EEE Now we recall that, by  \eqref{f:cc}, the set $\partial^ * \Omega _\flat  \cap (H _ T ^+ \setminus H _ T)$ is connected. 
   Hence for every $n \geq 1$ we can join $x_{1,n}$ to $x_{1,n+1}$  by a continuous arc $\gamma _{1,n}$ contained into $\partial ^* \Omega _\flat  \cap (H _ T ^+ \setminus H _ T)$. We can repeat the same  procedure for the second sequence, 
   constructing a family of continuous arcs $\gamma _{2,n}$  joining $x_{2,n}$ to $x_{2,n+1}$ for every $n \geq 1$. 
  
   We look at the boundaries of the balls of radius $r(s_0)$ whose centre moves along $\gamma _{1,n}$ and $\gamma_ {2, n}$. Clearly these balls tends to superpose in the limit as $n \to + \infty$, since $\|x_{1,n}-x_{2,n}\|$ decreases to $0$. 
   Moreover, we know from
  \eqref{f:separation2} that, during the continuous movement of their centre along 
  along $\gamma _{1,n}$ and  $\gamma _{2,n}$, 
the boundary of these balls cannot cross points of density $1$ for $\Omega ^ {\nu,  ns}$. 
Hence, for $n$ large,    
   $$
 B_\varepsilon (p) \cap \Omega ^ {\nu, ns} \sq  B_{r(s_0)} (x_{1,n}) \Delta B_{r(s_0)} (x_{2,n}) \,;$$
   hence, still for $n$ sufficiently large, 
$$|B_\varepsilon (p) \cap \Omega ^ {\nu, ns} | \leq | B_{r(s_0)} (x_{1,n}) \Delta B_{r(s_0)} (x_{2,n}) | < \frac{1}{4}  | B_\varepsilon (p)| \, , $$
against \eqref{f:alto}.

\smallskip In case $\overline x \in H _ T$, we proceed in the same way, except that  we cannot ensure any more that both  sequences $\{x_{1,n}\}$ and $\{x_{2,n}\}$ belong to the same halfspace $H _ T ^+$ or $H _ T ^-$. 
Thus, when we construct the continuous arcs $\gamma _{1,n}  $  and $\gamma _{ 2,n}$, they may belong indistinctly to $\partial ^* \Omega _\flat  \cap ( H ^ - _ T \setminus H_T)$ or to
$\partial ^* \Omega _\flat \cap ( H ^ + _ T \setminus H_T)$, but  this does not affect the validity of the proof since the contradiction follows as soon as $x_{1,n}$ and $x_{2,n}$ are close enough. \qed \EEE

\bigskip

    \begin{lemma}\label{bfm04} Under  the assumptions of Theorem \ref{t:serrin3},  given $\nu \in \mathbb S ^ { d-1}$, let $\Omega _ \flat$ be a  $h$-isolated open connected component  of $\Omega ^ {\nu, s}$. 
  There exists a constant $c_\flat >0$ such that
   \begin{equation}\label{f:ci} 
\int_    { \Om ^{\nu, s}} h _ x =c_\flat \qquad \forall x \in \overline{\partial^* \Om _\flat}\,.
   \end{equation} 
   
   \smallskip
Moreover, the constant is the same for any other open connected component   of $\Omega ^ {\nu, s}$  such that the closure of its essential boundary intersects 
$ \overline{\partial^* \Om_\flat}$.  
   \end{lemma} 
   
\proof   We argue in a similar way as in the proof of the previous lemma. 
Given  $
   x_1, x_2 \in \overline { \partial^* \Om_\flat} \cap (H_T^+\setminus H _ T) $,  
    by  \eqref{f:cc}, they can be joined by a continuous arc $\gamma$ contained into $\partial ^* \Omega _\flat \cap (H _ T ^+ \setminus H _ T)$.      
By  \eqref{f:separation2},  for $\mathcal L ^ 1$-a.e.\ $s \in (\sigma, + \infty)$, the boundary of the ball of radius $r(s)$ centred at any  point along $\gamma$
    cannot cross points of density $1$ for $\Omega ^ {\nu,  ns}$. 
   We deduce that, still for $\mathcal L ^ 1$-a.e.\ $s \in (\sigma, + \infty)$,   $B _ {r(s)} (x_1) \Delta B _ {r(s)} (x_2)$ cannot contain points of density $1$ for $\Omega ^ {\nu,  ns}$.    Therefore, by the equality \eqref{f:cake2} in Lemma \ref{l:cake}, we get
      $$ \int_ { \Om^{\nu, ns}  } |h _  {x_1} - h _  {x_2} | = 0    \,,$$  
and hence, using also the fact  that $\Omega$ is $h$-critical, 
$$\int _{\Om ^ {\nu , s }} h _ {x_1} = \int_ {\Om } h _ {x_1}   - \int_ { \Om^{\nu, ns}  } h _  {x_1}  = 
 \int_ {\Om } h _ {x_2}   - \int_ { \Om^{\nu, ns}  } h _  {x_2} =  \int_ {\Om ^ {\nu , s }} h _ {x_2}\,.$$ 

      By the arbitrariness of $x_1$, $x_2$, we infer that   there exists a   constant $c_\flat ^+>0$ such that
   $ \int_ { \Om^{\nu, s}  } h _ x  =c_\flat ^+$ for every $x \in \overline { \partial^* \Om_\flat} \cap( H^+_T \sm H_T)$. 
In the same way, we obtain that there exists a   constant $c_\flat ^->0$ such that
 $ \int_ { \Om^{\nu, s}  } h _ x  =c_\flat ^-$  for every $x \in \overline { \partial^* \Om_\flat} \cap( H^-_T  \sm H_T)$. 
 Since the two sets $\overline { \partial^* \Om _\flat } \cap H^{\pm}_T  $ 
  have common points on $H_T$, we conclude that $c_\flat ^ + = c_ \flat ^-$.
 The same argument proves also the last assertion of the lemma.  
   \qed

\bigskip \begin{lemma}\label{l:remove} 
Under  the assumptions of Theorem \ref{t:serrin3}, 
given $\nu \in \mathbb S ^ { d-1}$, let $\Omega _ \flat$ be a  $h$-isolated open connected component  of $\Omega ^ {\nu, s}$. 
  Then $\Omega _\flat$  is  and $h$-critical and not $h$-degenerate. The same assertions hold true for its complement $\Omega \setminus 
\Omega _\flat$,  unless it is of measure zero. 
\end{lemma}

\proof We know from  Lemma \ref{bfm03} that 
$$
 \inf_{x_1,x_2 \in \partial ^* \Om_\flat} \frac{\int_{\Om^{\nu, s}} | h_ {x_1}- h_{x_2}|}{\|x_1-x_2\|} >0.
 $$ 
But, since $\Omega_\flat$ is $h$-isolated,  we have 
$$ \int_{\Om^{\nu, s} }| h_ {x_1}- h_{x_2}|  =  \int_{\Om_\flat} | h_ {x_1}- h_{x_2}|  \qquad  \forall x_1,x_2 \in \partial ^* \Om_\flat $$
and hence $\Omega _ \flat$ is not $h$-degenerate. 

Let us prove that $\Omega _\flat$ is $h$-critical. For every $x \in \overline{\partial^* \Om_\flat}$, recalling  equality \eqref{f:ci} in Lemma \ref{bfm04}, we have
$$\int_    { \Om ^{\nu, s}} h _ x   = \sigma |\Omega ^ { \nu, s} | + 
 \int _{\sigma} ^ { + \infty} |\Omega^{\nu, s} \cap B _ { r (s) }(x)| \, ds
=  c _ \flat \,. 
$$
We infer that 
$$ \begin{array}{ll}\displaystyle \int_    { \Om _\flat } h _ x   & \displaystyle =  \sigma |\Omega _\flat| 
 + 
  \int _{\sigma} ^ { + \infty} |\Omega _\flat \cap B _ { r (s) }(x)| \, ds   
  \\  \noalign{\bigskip} 
  & \displaystyle =  \sigma |\Omega _\flat| 
 + 
  \int _{\sigma} ^ { + \infty} |\Omega ^ { \nu, s} \cap B _ { r (s) }(x)| \, ds
   =  \sigma |\Omega _\flat| 
 + c _ \flat - \sigma |\Omega ^ { \nu, s} |
    \,.
  \end{array}
$$ 
where the second equality follows from \eqref{f:separation}. 

\smallskip
Let us now consider the complement  $\Omega \setminus 
\Omega _\flat$.  Assume it is of positive measure, and hence that $\partial ^* ( \Omega \setminus 
\Omega _\flat )$ is not empty.

For every $x_1, x_ 2 \in \partial ^* (\Omega \setminus \Omega _\flat)$, by \eqref{f:separation2} and \eqref{f:separation},  it holds 
 
$$ \begin{array}{ll}  \displaystyle \int _{\Omega \setminus \Omega _ \flat } |h _ { x_1} - h _ {x_2} |  & \displaystyle = 
\int _{\sigma} ^ { + \infty}  | (\Omega \setminus \Omega _ \flat )  \cap ( B _  {r(s)}  ( x_1) \Delta  B _  {r(s)}  ( x_2) ) |  \, ds 
\\
\noalign{\medskip} 
&  \displaystyle =  
\int _{\sigma} ^ { + \infty}  | \Omega \cap ( B _ {r(s)}  ( x_1) \Delta  B _  {r(s)}  ( x_2) ) |  \, ds = 
\int _{\Omega  } |h _ { x_1} - h _ {x_2} | \,, 
\end{array}$$  
Thus
 $\Omega \setminus \Omega _\flat$ is not $h$-degenerate  since by assumption $\Omega$ is not $h$-degenerate.
 
In a similar way,  for every $x \in \partial ^* ( \Omega  \setminus \Omega _\flat)$, since $\Omega$ is $h$ critical, we have
$$\int_    { \Om } h _ x   = \sigma |\Omega  | + 
 \int _{\sigma} ^ { + \infty} |\Omega \cap B _ { r (s) }(x)| \, ds
=  c  \,. 
$$
We infer that
$$ \begin{array}{ll}  \displaystyle \int _{\Omega \setminus \Omega _ \flat } h _ {  x}   & \displaystyle = 
\sigma  | \Omega \setminus \Omega _ \flat | +  \int _{\sigma}  ^ { + \infty}  | (\Omega \setminus \Omega _ \flat )  \cap  B _ {r(s)}  ( x ) |  \, ds 
\\
\noalign{\medskip} 
&  \displaystyle =  \sigma |\Omega \setminus \Omega _\flat| + 
\int _{\sigma} ^ { + \infty}  | \Omega \cap  B _ {r(s)}  ( x)  |  \, ds = 
 \sigma |\Omega \setminus \Omega _\flat|  + c - \sigma |\Omega | \, ,
\end{array}$$  
where the second equality follows from 
\eqref{f:separation2} and \eqref{f:separation}.
 \qed

\bigskip

\medskip
{\it Proof of Claim 4b.}  First let us observe that the family $\mathcal F$ is at most countable. 
This follows from the facts that any open set of $\R ^d$ has at most countable connected components, 
and that, for two different directions $\nu_1$ and $\nu_ 2$, 
a connected component of  $\Om^{\nu _1, s}$ cannot intersect another connected component of  $\Om ^ { \nu _2, s} $ without being equal.  

\smallskip
We now prove Claim 4b. by contradiction. 
If the family $\mathcal F$ is not empty,
 it turns out to contain an element $\Omega  _ \sharp $ which is Steiner symmetric about $d$ hyperplanes 
with linearly independent normals $\nu _1, \dots, \nu _ d$ (for the detailed justification, see Claim 4b.\ in the proof of
\cite[Theorem 1]{BF}).

\smallskip Next we consider 
any other element $\Omega _ {\sharp \sharp}$  of $\mathcal F$ which is in $h$-contact with $\Omega _\sharp$ in the decomposition with respect to one among the directions $\nu_1 , \dots, \nu _d$,  say $\nu _ 1$. 
If $T_1$ is the stopping time for the parallel movement with normal $\nu _ 1$, 
there exist $p, p' \in \overline{\partial ^* \Om_\sharp} \setminus H_{T_1} $, symmetric about $H _{T_1}$,  such that
$$
\int_{\Om _{\sharp\sharp} } | h _ p - h _{p'} | = \int _{\sigma} ^{ + \infty}  \big |\big ( B_{r(s)} (p)\Delta B_{r(s)} (p')\big ) \cap  \Om_{\sharp\sharp}  \big | \, ds >0.
$$ 
We infer that
$$\mathcal L ^ 1 \big ( \big \{ s \in(\sigma, + \infty)\ :\ \big |\big ( B_{r(s)} (p)\Delta B_{r(s)} (p')\big ) \cap  \Om_{\sharp\sharp}  \big | >0 \big \} \big ) >0 \,,$$ 
and hence, since we are assuming that $\Om_{\sharp\sharp} $ is   Steiner symmetric with respect to $H_{T_1}$, 
$$\mathcal L ^ 1 \big ( \big \{ s \in (\sigma, + \infty)\ :\  \partial B _ {r(s)} (p) \cap   \Om^ { (1)}_{\sharp\sharp}  \neq \emptyset  \big \} \big ) >0 \,.$$

Recalling \eqref{f:diff-balls}, this implies that 
$\Omega _ {\sharp \sharp}$  is itself Steiner symmetric about the same hyperplanes as $\Omega _\sharp$ is. 
Then, Lemma 18 in \cite{BF}  \EEE implies that the set $\Om_\sharp \cup \Om_{\sharp\sharp}$ is connected, yielding a contradiction.

\medskip
{\it Proof of Claim 4c.} 
 We start the procedure by choosing a direction $\nu \in \mathbb S ^ { d-1}$. 
By Claim 4b., we can pick a $h$-isolated open connected component of $\Omega ^ {\nu, s}$, which  by Claim 4a. turns out to be a ball $B_1$ of radius $R_1> \eta  /2$. We remove this ball from $\Omega$. 
By Claim 4a.,  
we are left with a set $\Omega'$ which is still $h$-critical and  not $h$-degenerate  (unless it has measure zero).  So we can restart the process
with  $\Omega '$ in place of $\Omega$. 
Again, by Claim 4b., we can pick a $h$-isolated open connected component of $(\Omega  ') ^ {\nu, s} $, which  by Claim 4a. turns out to be a ball $B_2$ of radius $R_2> \eta  /2$. We remove this ball from $ \Omega'$.  We observe that, 
by \eqref{f:separation}, we have  
\begin{equation}\label{f:bound-dist} {\rm dist} (B _1, B _ 2) \geq  \sup_ { s \in(\sigma, + \infty) } r(s) = 
r (\sigma)  \,.
\end{equation} 
Then, for every $p _ 1 \in \partial B _ 1$ and $p _ 2 \in \partial B _ 2$,  we have
$$\begin{array}{ll} \displaystyle \int _ {\sigma} ^ { + \infty} |B _ 1 \cap B _ { r(s)} ( p _1)| \, ds& \displaystyle =
 \int _ {\sigma}  ^ { + \infty} |(B _ 1 \cup B _ 2)  \cap B _ { r(s)} ( p _1) | \, ds  = \\ \noalign{\bigskip}
 & = \displaystyle
 \int _ {\sigma}  ^ { + \infty} |(B _ 1 \cup B _ 2)  \cap B _ { r(s)} ( p _2) | \, ds  = 
   \int _ {\sigma}  ^ { + \infty} |B _ 2 \cap B _ { r(s)} ( p _2) | \, ds \,. \end{array}$$ 
where the first and third equalities follow from \eqref{f:bound-dist}, while the second one is consequence of the $h$-criticality of $\Omega$ and of $\Omega \setminus ( B _ 1 \cup B _ 2)$. 
Taking into account that,  for every fixed $s$, the map $R \mapsto  | B _ R \cap B _ { r(s)} (p) |$, with $p \in \partial B _ R$, is strictly increasing, we deduce that $R_ 1 = R _ 2$. 

Since $\Omega$ has finite measure, we can repeat this process a finite number of times, until when
we are left with a set of measure zero.   \qed

  \bigskip

\bibliographystyle{mybst}

\end{document}